\documentclass[12pt,a4paper]{amsart}
 \usepackage[T1]{fontenc}
\usepackage[top=35mm, bottom=35mm, left=30mm, right=30mm]{geometry}    
\usepackage{amssymb,amsmath} 
\usepackage{mathtools}
\usepackage{verbatim,enumerate}

\usepackage[looser]{newpxtext}
\usepackage[smallerops]{newpxmath}

\usepackage{fancyhdr}
\makeatletter
\def\@@and{\MakeLowercase{and}}
\makeatother

\usepackage[hidelinks,pagebackref]{hyperref}

\renewcommand*{\backref}[1]{}
\renewcommand*{\backrefalt}[4]{%
\ifcase #1 No citations.%
\or (Cited on page #2)%
\else (Cited on pages #2)%
\fi%
}

\theoremstyle{definition}
\newtheorem{defn}{Definition}[section]
\newtheorem{exam}[defn]{Example}

\newtheorem{rem}[defn]{Remark}
\newtheorem{Ques}[defn]{Question}

\theoremstyle{plain}
\newtheorem{thm}[defn]{Theorem}
\newtheorem{lem}[defn]{Lemma}
\newtheorem{prop}[defn]{Proposition}
\newtheorem{coro}[defn]{Corollary}

\title[W\MakeLowercase{eak disjointness of hypercyclic operators}]
{W\MakeLowercase{eak disjointness of hypercyclic operators}}

\author[J. L\MakeLowercase{i}]{J\MakeLowercase{ian} L\MakeLowercase{i}}
\address[J. Li]{Institute for  Mathematical Sciences and Artificial Intelligence \& Department of Mathematics,
	Shantou University, Shantou, 515821, Guangdong, China}
\email{lijian09@mail.ustc.edu.cn}
\urladdr{https://orcid.org/0000-0002-8724-3050}

\author[Q. L\MakeLowercase{iao}]{Q\MakeLowercase{ijing} L\MakeLowercase{iao}}
\address[Q. Liao]{Department of Mathematics,
	Shantou University, Shantou, 515821, Guangdong, China}
\email{liaoqijing1@outlook.com}

\author[Y. R\MakeLowercase{uan}]{Y\MakeLowercase{onghang} R\MakeLowercase{uan}}
\address[Y. Ruan]{Department of Mathematics,
	Shantou University, Shantou, 515821, Guangdong, China}
\email{yonghangruan@outlook.com}

	

\date{\today}

\allowdisplaybreaks[3] 
 
\begin{document}

\begin{abstract}
We study the weak disjointness of hypercyclic operators to refine the classifications of hypercyclic operators.
We establish an analogue of the Weiss-Akin-Glasner Theorem from topological dynamics within the framework of linear dynamics, which gives a characterization of the weak disjointness of each class of mixing operators with respect to a given Furstenberg family.
The key ingredient is an analogue of the Weiss-Akin-Glasner Lemma from topological dynamics, which gives a characterization of subsets of non-negative integers which can be realized by the return time sets of mixing  operators with respect to a given  Furstenberg family.
We also provide several examples to distinguish some classes of hypercyclic operators and end with the characterization of the mixing property with respect to Furstenberg families of direct sums of backward shifts on Fr\'echet sequence spaces. 
\end{abstract}

\maketitle

\section{Introduction}

In 1967, Furstenberg introduced the concept of disjointness in topological dynamics and ergodic theory in the seminal paper \cite{F67}, and characterized certain classes of systems by their disjointness relations. 
Two compact dynamical systems $(X,T)$ and $(Y,S)$ (continuous self-maps on compact metric spaces) are said to be disjoint if $X\times Y$ is the only closed $T\times S$-invariant subset whose coordinate projections are surjective.
It is easy to see that if two transitive compact dynamical systems $(X,T)$ and $(Y,S)$ are disjoint, then the product system $(X\times Y,T\times S)$ is transitive.
Later in 1972, Peleg introduced the concept of weak disjointness in topological dynamics in~\cite{P72}, which, as the name suggests, is a condition weaker than disjointness. 
Two compact dynamical systems $(X,T)$ and $(Y,S)$ are said to be weakly disjoint provided that the product system $(X\times Y,T\times S)$ is  transitive.

In 2000, Weiss characterized in \cite{W00} the compact dynamical systems which are weakly disjoint from all weakly mixing compact dynamical systems, as being exactly the syndetically transitive ones.
Let $\mathcal{F}$ be a Furstenberg family (a hereditary upward collection of subsets of non-negative integers). 
A topological dynamical system $(X,T)$ is said to be $\mathcal{F}$-transitive provided that $N(U,V):=\{n\in\mathbb{N}_0: T^nU\cap V\neq\emptyset\}\in\mathcal{F}$ for any non-empty open subsets $U,V$ of $X$; 
and $\mathcal{F}$-mixing provided that the product system $(X\times X,T\times T)$ is $\mathcal{F}$-transitive. 
In 2001, Akin and Glasner studied in~\cite{AG01} the $\mathcal{F}$-transitivity of compact dynamical systems and
provided a method to characterize the compact dynamical systems which are weakly disjoint from a prescribed class of compact dynamical systems, generalizing Weiss' result. 
The key ingredient is the realization of every set in a shift-invariant thick Furstenberg family $\mathcal{F}$ by the return time sets of $\mathcal{F}$-transitive compact dynamical systems.
The corresponding result of the characterization of compact dynamical systems is known as the Weiss-Akin-Glasner Theorem, and the realization result is called the Weiss-Akin-Glasner Lemma, see e.g.\@ \cite{HY02} and \cite{HY04}. 
It turns out that the notion of weak disjointness plays a significant role in the classification of transitive compact dynamical systems, see e.g.\@ \cite{AG01,G04,HY02}.

Furstenberg showed in \cite{F67} that the return time sets of any weakly mixing topological dynamical system form a filter base. 
By this result it is easy to see that a topological dynamical system $(X,T)$ is weakly mixing if and only if $N(U,V)$ is thick for any non-empty open subsets $U$ and $V$ of $ X$, and $(X,T)$ is $\mathcal{F}$-mixing if and only if it is $\mathcal{F}$-transitive and weakly mixing.
This implies that for a thick family $\mathcal{F}$, the $\mathcal{F}$-transitivity is equivalent to the $\mathcal{F}$-mixing. 
Then the $\mathcal{F}$-transitive system in the Weiss-Akin-Glasner Lemma is in fact $\mathcal{F}$-mixing. 
On the other hand, the family generated by the return time sets of a non-trivial $\mathcal{F}$-mixing topological dynamical system must be a shift-invariant thick family.  
Therefore, the Weiss-Akin-Glasner Lemma gives a characterization of subsets of non-negative integers which can be realized by the return time sets of $\mathcal{F}$-mixing compact dynamical systems (see Lemma~\ref{lem:WAG-revised}), and the Weiss-Akin-Glasner Theorem characterizes the compact dynamical systems which are weakly disjoint from $\mathcal{F}$-mixing compact dynamical systems (see Theorem~\ref{thm:WAG-revised}).

Linear dynamics has become an active area of research in the last decades, see the monographs \cite{B09} and \cite{G11} for the introduction of dynamics of linear operators.
In particular, the study of hypercyclic operators constitutes a subject of considerable interest. 
A linear dynamical system (a continuous linear operator on an infinite-dimensional Fr\'echet space over $\mathbb{R}$ or $\mathbb{C}$) is said to be hypercyclic provided that the orbit of some vector is dense. 
Although linear dynamical systems lack compactness compared to compact dynamical systems, by different methods and techniques many conclusions in linear dynamical systems are parallel to those in compact dynamical systems. 

As every linear dynamical system has a fixed point, two linear dynamical systems can never be disjoint in the sense of Furstenberg. 
In \cite{B07} and \cite{BP07}  the authors introduced the notion of disjoint hypercyclicity, which is weaker than disjointness. 
In \cite{B07,BP07,BMPP19,C24,S95}, the authors studied the hypercyclicity of direct sums, which is in fact the weak disjointness in linear dynamics. 
In 1992, Herrero asked in \cite{H92} whether $T\oplus T$ is hypercyclic or not provided that $T$ is hypercyclic. 
In 1995, Salas gave an affirmative answer in \cite{S95} for the case where $T$ is a weighted backward shift on $\ell^2(\mathbb{N}_0)$ or $\ell^2(\mathbb{Z})$, and also provided two hypercyclic operators $T,S$ on $\ell^2(\mathbb{Z})$ such that $T\oplus S$ is not hypercyclic. 
De la Rosa and Read \cite{DR09}
constructed a hypercyclic operator on a suitable Banach space such that the direct sum is not hypercyclic, while Bayart and Matheron 
\cite{BM07} showed that operators of this type can be found on many classical Banach spaces like $c_0(\mathbb{N}_0)$, $\ell^p(\mathbb{N}_0)$, $1\leq p<\infty$, etc.
To provide more hypercyclic operators without hypercyclic direct sums, it is natural to ask: which operators are weakly disjoint from a given class of hypercyclic operators? 

In \cite{BMPP19}, B\`es et al.\@ studied the $\mathcal{F}$-transitivity of operators. 
They formulated a criterion for an operator to be $\mathcal{F}$-transitive and characterized the $\mathcal{F}$-transitivity of weighted backward shifts. 
In \cite{C24}, Cardeccia studied the hypercyclicity of the direct sum of an operator and its iterates.  
Cardeccia also proved that the operators which are weakly disjoint from all weakly mixing operators are exactly the syndetically transitive ones, and that the operators which are weakly disjoint from all syndetically transitive operators are exactly the piecewise syndetically transitive ones. 

The aim of this paper is to study the weak disjointness of hypercyclic operators, advancing the classifications of hypercyclic operators.
Motivated by the Weiss-Akin-Glasner Theorem from topological dynamics and Cardeccia's results for operators, we characterize the weak disjointness of more classes of operators. 
We first present some preliminaries in Section 2.
In Section 3, we state some basic properties of weak disjointness in linear dynamics, and give sufficient conditions for $\mathcal{F}$-mixing operators (Theorem \ref{F-mixing1}) and mild mixing operators (Proposition \ref{mild-ip}). 
We also illustrate some examples to show that a weakly mixing operator may not be weakly disjoint from its iterates.
In Section 4, as an analogue of the Weiss-Akin-Glasner Theorem from topological dynamics,
we formulate the corresponding result (Theorem \ref{linearWAG}) in linear dynamics, which gives a characterization of the weak disjointness of $\mathcal{F}$-mixing operators.
The key ingredient of the proof is an analogue of the Weiss-Akin-Glasner Lemma from topological dynamics, see Lemma \ref{L-WAG}, which gives a characterization of subsets of non-negative integers which can be realized by the return time sets of $\mathcal{F}$-mixing operators. 
We also provide several examples of operators with some prescribed properties. 
In Section 5, we characterize the $\mathcal F$-mixing property of direct sums of backward shifts on Fr\'echet sequence spaces in terms of the F-norms of the basis. 

\section{Preliminaries}

In this section we recall some basic notions and results that will be used consistently throughout the paper.

\subsection{Subsets of non-negative integers} 
By $\mathbb{N}$, $\mathbb{N}_0$ and $\mathbb{Z}$, we denote respectively the set of positive integers, the set of non-negative integers, and the set of integers. 

Before proceeding further, we briefly recall the relevant notations concerning a Furstenberg family (detailed discussions can be found in \cite{A97} and \cite{F81}). 
Let $\mathcal{F}$ be a collection of subsets of $\mathbb{N}_0$. 
We say that $\mathcal{F}$ is a \emph{Furstenberg family} (or \emph{family} for short) provided that it is hereditary upward, that is, for any $A\in \mathcal{F}$ and $B\subset\mathbb{N}_0$ with $A\subset B$  we have $B\in \mathcal{F}$. 
Denote the family of all infinite subsets of $\mathbb{N}_0$ by $\mathcal{F}_{\mathrm{inf}}$ and the family of all cofinite subsets of $\mathbb{N}_0$ by $\mathcal{F}_{\mathrm{cf}}$. 
For a family $\mathcal{F}$, we say that a collection $\mathcal{F}'$ of subsets of $\mathbb{N}_0$ is a \emph{subfamily} of $\mathcal{F}$ if $\mathcal{F}'$ is a family with $\mathcal{F}'\subset\mathcal{F}$. 

\textbf{Throughout this paper, for each Furstenberg family $\mathcal{F}$, we adopt the conventions that $\emptyset\notin\mathcal{F}$ and $\mathcal{F}\neq \emptyset$. }

For a subset $F$ of $\mathbb{Z}$ and $i\in\mathbb{Z}$, let
$F\pm i=\{m\pm i\colon m\in F\}$. 
We say that a family $\mathcal{F}$ is 
\begin{itemize}
\item \emph{shift$_+$-invariant} if $F+i\in\mathcal{F}$ for each $i\in\mathbb{N}_0$ and $ F\in\mathcal{F}$; 
\item \emph{shift$_-$-invariant} if $(F-i)\cap \mathbb{N}_0\in\mathcal{F}$ for each $i\in\mathbb{N}_0$ and $ F\in\mathcal{F}$; 
\item \emph{shift-invariant} if $(F+i)\cap\mathbb{N}_0 \in\mathcal{F}$ for each $i\in\mathbb{Z}$ and $F\in\mathcal{F}$. 
\end{itemize}

\begin{defn}
Let $\mathcal{F}$ be a family. 
The \emph{dual family} $\mathcal{F}^*$ is given by: 
\[
\mathcal{F}^*=\{A\subset \mathbb{N}_0: A\cap B\neq\emptyset \text{ for all }B\in\mathcal{F}\}. 
\]
We say that $\mathcal{F}$ is a \emph{full family} if $F\cap F'$ is infinite for any $F\in \mathcal{F}$ and $F'\in\mathcal{F}^*$; 
$\mathcal{F}$ is a \emph{filter} if $F_1\cap F_2\in \mathcal{F}$ for any $F_1, F_2\in \mathcal{F}$. 
\end{defn}

\begin{rem}
The following statements can be checked easily.
    \begin{enumerate}
        \item A family $\mathcal{F}$ is full if and only if $F\cap F'\in \mathcal{F}$ for any $F\in \mathcal{F}$ and $F'\in \mathcal{F}_{\mathrm{cf}}$;
        \item each shift-invariant family is full; 
        \item for each full family $\mathcal{F}$, we must have $\mathcal{F}_{\mathrm{cf}}\subset \mathcal{F}$; 
        \item a family $\mathcal{F}$ is full if and only if $\mathcal{F}^*$ is full. 
    \end{enumerate}
\end{rem}

We recall some standard families: 
\begin{itemize}
    \item $\mathcal{F}_{\mathrm{t}}$: the family of subsets of $\mathbb{N}_0$ containing intervals of arbitrary length, which are said to be \emph{thick};  
    \item $\mathcal{F}_{\mathrm{s}}$: the family of subsets of $\mathbb{N}_0$ with bounded gaps, which are said to be \emph{syndetic}; 
    \item $\mathcal{F}_{\mathrm{ps}}$: the family of subsets of $\mathbb{N}_0$ with the form $F_1\cap F_2$, where $F_1\in\mathcal{F}_{\mathrm{s}}$ and $F_2\in\mathcal{F}_{\mathrm{t}}$, and these are said to be \emph{piecewise syndetic};  
    \item $\mathcal{F}_{\mathrm{ts}}$: the family of subsets $A$ of $\mathbb{N}_0$ such that for each $k\in\mathbb{N}_0$, there exists $F\in\mathcal{F}_{\mathrm{s}}$ with $\bigcup_{i=0}^k (F+i)\subset A$, and these are said to be \emph{thickly syndetic}. 
\end{itemize}
It is easy to see that $\mathcal{F}_{\mathrm{s}}^*=\mathcal{F}_{\mathrm{t}}$ and $ \mathcal{F}_{\mathrm{ps}}^*=\mathcal{F}_{\mathrm{ts}}$. 

We turn to a property of families (see \cite{A97,AG01}) which plays a significant role in our later discussion.  

\begin{defn}
Let $\mathcal{F}$ be a family. 
We say that $\mathcal{F}$ is a \emph{thick family} if 
for any $F\in\mathcal{F}$ and $k\in\mathbb{N}_0$, 
\[
    \bigcap_{i=0}^k(F-i)\cap\mathbb{N}_0\in\mathcal{F}.
\]
\end{defn}

\begin{rem}\label{thick}
We have the following easy facts:
\begin{enumerate}
\item each thick family is shift$_-$-invariant; 
\item each thick family is a subfamily of $\mathcal{F}_{\mathrm{t}}$; 
\item each shift-invariant filter is a thick family;   
\item a family $\mathcal{F}$ is a shift-invariant thick family if and only if for each $F\in\mathcal{F}$ and any $k_1, k_2\in\mathbb{N}_0$, 
       \[
           \bigl\{n\in\mathbb{N}_0: \{n-k_1,n-k_1+1,\dotsc,n+k_2\}\subset F\bigr\}\in\mathcal{F};
       \]
\item both $\mathcal{F}_{\mathrm{t}}$ and $\mathcal{F}_{\mathrm{ts}}$ are shift-invariant thick families. 
\end{enumerate}
\end{rem}

Following \cite{BMPP19}, for a given family $\mathcal{F}$, 
we denote 
\[
    \widetilde{\mathcal{F}}=\biggl\{F\in\mathcal{F}: \text{for any } k\in\mathbb{N}_0, \bigcap_{i=-k}^k(F+i)\cap\mathbb{N}_0\in\mathcal{F} \biggr\}. 
\]
Then $\widetilde{\mathcal{F}}$ is the maximal shift-invariant thick subfamily of $\mathcal{F}$ provided that it is not empty. 
Note that if $\mathcal{F}$  is full, then $\mathcal{F}_{\mathrm{cf}}\subset \widetilde{\mathcal{F}}$.

\subsection{Topological dynamics} 
For a \emph{topological dynamical system}, 
we mean a pair $(X,T)$, where  $X$ is a metric space and $T$ is a continuous map of $X$ into $X$.
If in addition $X$ is a compact metric space, then we say that $(X,T)$ is a \emph{compact dynamical system}. 
See \cite{F67, F81} for details of compact dynamics.

A topological dynamical system $(X,T)$ is said to be \emph{topologically transitive} (or \emph{transitive} for short) if for any pair $U, V$ of non-empty open subsets of $X$, the return time set 
\[
    N_T(U,V)=N(U,V):=\{n\in\mathbb{N}_0: T^nU\cap V\neq\emptyset \}
\]
is non-empty. 
When the operator is fixed and there is no risk of ambiguity, we shall write $N(U,V)$ in place of $N_T(U,V)$.
A topological dynamical system $(X,T)$ is said to be 
\emph{weakly mixing} if the product system $(X\times X,T\times T)$ is transitive, and \emph{strongly mixing} if for any pair $U, V$ of non-empty open subsets of $X$, $N(U,V)\in\mathcal{F}_{\mathrm{cf}}$.

Given a family $\mathcal{F}$, we say that a topological dynamical system $(X,T)$ is 
\begin{itemize}
    \item \emph{$\mathcal{F}$-transitive} if for any non-empty open subsets $U,V$ of $X$, $N(U,V)\in\mathcal{F}$; 
    \item \emph{$\mathcal{F}$-mixing} if the product system $(X\times X,T\times T)$ is $\mathcal{F}$-transitive. 
\end{itemize}
For each family $\mathcal{F}$, since $\emptyset\notin\mathcal{F}$, the $\mathcal{F}$-transitivity is stronger than topological transitivity. 
It can be verified that a topological dynamical system is transitive if and only if it is $\mathcal{F}_{\mathrm{inf}}$-transitive. 
A topological dynamical system is said to be \emph{topologically ergodic} if it is $\mathcal{F}_{\mathrm{s}}$-transitive. 

We have the following characterization of weak mixing of topological dynamical systems, which is essentially contained in the proof of \cite[Proposition II.3]{F67}.
\begin{prop}\label{F-intersection}
Let $(X,T)$ be a topological dynamical system. 
Then $(X,T)$ is weakly mixing if and only if it is $\mathcal{F}_{\mathrm{t}}$-transitive and if and only if 
\[
    \mathcal{N}_T=\bigl\{A\subset\mathbb{N}_0\colon N(U,V)\subset A\text{ for some non-empty open subsets }U,V\text{ of }X\bigr\}
\]
is a filter. 
\end{prop}

By Proposition~\ref{F-intersection},
it is easy to see that a topological dynamical system $(X,T)$ is $\mathcal{F}$-mixing if and only if it is $\mathcal{F}$-transitive and weakly mixing. 
We have the following folklore result on the return time sets of $\mathcal{F}$-mixing  systems, see e.g.\@
\cite[Lemma 2.3]{BMPP19} in the setting of linear dynamics.

\begin{prop}\label{return-time-F-mixing}
Let $(X,T)$ be a topological dynamical system where $X$ is not a singleton and $\mathcal{F}$ be a family.
Then $(X,T)$ is $\mathcal{F}$-mixing if and only if 
\[
    \mathcal{N}_T=\bigl\{A\subset\mathbb{N}_0\colon N(U,V)\subset A\text{ for some non-empty open subsets }U,V\text{ of }X\bigr\}
\]
is a shift-invariant thick subfamily of $\mathcal{F}$. 
\end{prop}

\begin{proof}
On the one hand, suppose that $T$ is $\mathcal{F}$-mixing. 
It is clear that $\mathcal{N}_T$ is a subfamily of $\mathcal{F}$. 
By Proposition~\ref{F-intersection}, $\mathcal{N}_T$ is a filter. 
By Remark~\ref{thick}, it suffices to show that $\mathcal{N}_T$ is a shift-invariant family.
Let $U$ and $V$ be non-empty open subsets of $X$ and $i\in\mathbb{N}_0$. 
Since $T$ is transitive, it has a dense range and then $T^{-i}U\neq\emptyset$ and $T^{-i}V\neq\emptyset$. 
Since $T$ is weakly mixing and $X$ is not a singleton, we have that $X$ has no isolated points. 
Then there are non-empty open subsets $W_1\subset T^{-i}U, W_2\subset V$ with $N(W_1,W_2)\subset [i,+\infty)$. 
Note that $N(T^{-i}U,V)\cap [i,+\infty)\subset N(U,V)+i$. 
So we have 
\[
    N(W_1,W_2)\subset N(T^{-i}U,V)\cap [i,+\infty)\subset N(U,V)+i, 
\]
and 
\[
    N(U,T^{-i}V)\subset (N(U,V)-i)\cap\mathbb{N}_0. 
\]
By the fact that $\mathcal{N}_T$ is a filter, 
$\mathcal{N}_T$ is a shift-invariant thick family.

On the other hand, suppose that $\mathcal{N}_T$ is a shift-invariant thick subfamily of $\mathcal{F}$. 
Then we have $\mathcal{N}_T\subset\mathcal{F}_{\mathrm{t}}$. 
By Proposition \ref{F-intersection}, $T$ is weakly mixing. 
Since $T$ is $\mathcal{F}$-transitive, we have that $T$ is $\mathcal{F}$-mixing. 
\end{proof}

\begin{rem}
Let $\mathcal{F}_1$ and $\mathcal{F}_2$ be two families.
If there exists an $\mathcal{F}_1$-mixing $(X,T)$ which is not $\mathcal{F}_2$-transitive, then by Proposition \ref{return-time-F-mixing}, 
there exists a shift-invariant thick subfamily $\mathcal{F}_0$ of $\mathcal{F}_1$ such that $\mathcal{F}_0\setminus\mathcal{F}_2\neq\emptyset$.
\end{rem}

\begin{defn}
Two topological dynamical systems $(X,T)$ and $(Y,S)$ are said to be \emph{weakly disjoint} provided that their product $(X\times Y,T\times S)$ is topologically transitive. 
\end{defn}

In 2000, Weiss found in \cite{W00} that the compact dynamical systems which are weakly disjoint from all weakly mixing compact dynamical systems, are exactly the topologically ergodic ones. Later in 2001, Akin and Glasner generalized Weiss' result, to more classes of compact dynamical systems in \cite{AG01}, and the generalized result is known as the Weiss-Akin-Glasner Theorem, see e.g.\@ \cite{HY02} and \cite{HY04}. 
In this theorem, they characterized the weak disjointness of $\mathcal{F}$-transitive compact dynamical systems where $\mathcal{F}$ is a shift-invariant thick family, 
and connected the study of topological dynamics and combinatorial number theory. 

\begin{thm}[Weiss-Akin-Glasner Theorem, {\cite[Theorem 4.15]{AG01}}]\label{WAG-thm}
Let $(X,T)$ be a compact dynamical system and $\mathcal{F}$ be a shift-invariant thick family.
Then $(X,T)$ is $\mathcal{F}^*$-transitive if and only if  
$(X,T)$ is weakly disjoint from all $\mathcal{F}$-transitive compact dynamical systems. 
\end{thm}

The key ingredient of the proof of the Weiss-Akin-Glasner Theorem is the following lemma. 

\begin{lem}[Weiss-Akin-Glasner Lemma, {\cite[Lemma 4.16]{AG01}}]\label{WAG-lem}
Suppose that $\mathcal{F}$ is a shift-invariant thick family and $F\in\mathcal{F}$. 
Then there exists an $\mathcal{F}$-transitive compact dynamical system $(X,T)$, and a non-empty open subset $U$ of $X$ such that 
     \[
         N(U,U)=F\cup\{0\}.
     \]
\end{lem}

\begin{rem}
As $\mathcal{F}$ is a shift-invariant thick family, by Remark~\ref{thick} 
the $\mathcal{F}$-transitive compact dynamical system $(X,T)$ in the Weiss-Akin-Glasner Lemma is in fact $\mathcal{F}$-mixing. 

By Proposition~\ref{return-time-F-mixing},
the family generated by the return time sets of an $\mathcal{F}$-mixing system whose underlying space is not a singleton must be a shift-invariant thick family. 
So the condition that $\mathcal{F}$ is a shift-invariant thick family 
in the Weiss-Akin-Glasner Lemma is necessary for all sets in $\mathcal{F}$ to be realized by the return time sets of $\mathcal{F}$-mixing systems. 
\end{rem}

Recall that $\widetilde{\mathcal{F}}$ is the maximal shift-invariant thick subfamily of a full $\mathcal{F}$. Then we can restate the Weiss-Akin-Glasner Lemma as follows:
\begin{lem} \label{lem:WAG-revised}
Suppose that $\mathcal{F}$ is a full family and $F\subset\mathbb{N}_0$.
Then $F\in \widetilde{\mathcal{F}}$ if and only if there exists an $\mathcal{F}$-mixing compact dynamical system $(X,T)$, and a non-empty open subset $U$ of $X$ such that 
     \[
         N(U,U)=F\cup\{0\}.
     \]
\end{lem}
We also restate the Weiss-Akin-Glasner Theorem as follows:

\begin{thm}\label{thm:WAG-revised}
Let $(X,T)$ be a compact dynamical system and $\mathcal{F}$ be a full family.
Then $(X,T)$ is $\bigl(\widetilde{\mathcal{F}}\bigr)^*$-transitive if and only if  
$(X,T)$ is weakly disjoint from all $\mathcal{F}$-mixing compact dynamical systems. 
\end{thm} 

The concept of mild mixing was originally introduced by Furstenberg and Weiss in the context of ergodic theory \cite{FW06}, and was later introduced independently in topological dynamics by Glasner \cite{G04} and by Huang and Ye \cite{HY04}, inspired by its ergodic-theoretic counterpart.

\begin{defn}
    A compact dynamical system $(X,T)$ is said to be \emph{mildly mixing} if it is weakly disjoint from all transitive compact dynamical systems. 
\end{defn}

The return time sets of mild mixing systems are related to IP-sets. 
For a sequence $(x_n)_{n\in\mathbb{N}}$ in $\mathbb{N}$, 
the \emph{finite sum set} of $(x_n)_{n\in\mathbb{N}}$ is given by 
\[
\operatorname{FS}((x_n)_{n\in\mathbb{N}})=\{x_{i_1}+x_{i_2}+\dotsb+x_{i_n} \colon i_1<i_2<\dotsb <i_n, n\in \mathbb{N}\}.\]
A set $A\subset \mathbb{N}_0$ is called an \emph{IP-set} if there is a sequence $(x_n)_{n\in\mathbb{N}}$ in $\mathbb{N}_0$ with infinitely many non-zero terms such that $\operatorname{FS}((x_n)_{n\in\mathbb{N}})\subset A$. 
Denote by $\mathcal{F}_{\mathrm{ip}}$ the family of all IP-sets, and by $\Delta_{\mathrm{ip}}$ the family of the sets of the form $(F-F)\cap\mathbb{N}_0$, where $F\in\mathcal{F}_{\mathrm{ip}}$ and $F-F=\{m-n:m\neq n\in F\}$.
We denote the dual family of $\Delta_{\mathrm{ip}}$ by $\Delta^*_{\mathrm{ip}}$. For a topological dynamical system $(X,T), x\in X$ and a subset $U$ of $X$, we denote 
\[
N_T(x,U)=N(x,U):=\{n\in\mathbb{N}_0: T^nx\in U\}. 
\]
Recall that a point $x\in X$ is \emph{recurrent} 
if for every neighborhood $U$ of $x$, $N(x,U)$ is infinite.
In 1981, Furstenberg proved the following result for compact dynamical systems, which is in fact true for topological dynamical systems. 

\begin{thm}[{\cite[Theorem 2.17]{F81}}]\label{NxW-ip}
Let $(X, T)$ be a topological dynamical system. If a point $x\in X$ is recurrent, then for any neighborhood $U$ of $x$, we have $N(x, U)\in\mathcal{F}_{\mathrm{ip}}$. 
\end{thm}

The following characterization of mild mixing was obtained by 
Glasner \cite{G04}, and by Huang and Ye \cite{HY04}, independently. 

\begin{thm}[{\cite[Theorem 10.7]{G04} or \cite[Theorem 6.6]{HY04}}]\label{cpt-mm-ip}
A compact dynamical system is mildly mixing if and only if it is $\Delta^*_{\mathrm{ip}}$-transitive. 
\end{thm}

\subsection{Linear dynamics}
In recent years, there has been growing interest in the study of linear dynamical systems, see \cite{B09} and \cite{G11} for further information. 

Let $X$ be an infinite-dimensional separable Fr\'echet space over $\mathbb{R}$ or $\mathbb{C}$. 
For an \emph{operator} (i.e.\@ a continuous linear self-map) $T$ on $X$, the pair $(X,T)$ is called a \emph{linear dynamical system}. 
Usually we simply call the \emph{operator} $T$ or $T\colon X\to X$ a linear dynamical system. 
If in addition $T$ is an invertible map such that $T^{-1}$ is also an operator on $X$, then we say that $(X,T)$ (or $T$ for short) is \emph{invertible}. 
In this paper, for a Fr\'echet space $X$, the symbol $\|\cdot\|_X=\|\cdot\|$ denotes an F-norm defining the topology of $X$. 
If in particular $X$ is a Banach space, then the symbol $\|\cdot\|_X=\|\cdot\|$ denotes a norm defining the topology of $X$. 

While infinite-dimensional Fr\'echet spaces are lacking of compactness, many results in linear dynamical systems remain similar to their counterparts from compact dynamical systems, despite being obtained by different methods. 

\begin{defn}
An operator $T\colon X\to X$ is said to be \emph{hypercyclic} provided that there is some $x\in X$ whose orbit under $T$ is dense in $X$. 
Such an $x$ is called a \emph{hypercyclic vector} for $T$. We denote the set of hypercyclic vectors for $T$ by $HC(T)$.
\end{defn}

\begin{rem}\label{Birkhoff}
By the Birkhoff transitivity theorem (see e.g.\@ \cite[Theorem 1.16]{G11}) we know that for a separable infinite-dimensional Fr\'echet space $X$, an operator $T$ on $X$ is hypercyclic if and only if it is topologically transitive and if and only if $HC(T)$ is a dense $G_\delta$ subset of $X$.
\end{rem}

A well-known result to imply the weak mixing of a given operator is the so-called Hypercyclicity Criterion. 

\begin{thm}[Hypercyclicity Criterion, see e.g.\ {\cite[Theorem 3.12]{G11}}]\label{HC}
Let $T\colon X\to X$ be an operator. 
If there exist dense subsets $X_0$, $Y_0\subset X$, a sequence $(n_k)_{k\in\mathbb{N}}$ in $\mathbb{N}$ and a sequence of maps $S_{n_k}\colon Y_0\to X$ such that for any $x_0\in X_0$, $y_0\in Y_0$, 
   \begin{enumerate}
       \item $\displaystyle\lim_{k\to\infty}T^{n_k}x_0=0$; 
       \item $\displaystyle\lim_{k\to\infty}S_{n_k}y_0=0$; 
       \item $\displaystyle\lim_{k\to\infty}T^{n_k}S_{n_k}y_0=y_0, $
   \end{enumerate}
then $T$ is weakly mixing. 
\end{thm}

A remarkable feature of the Hypercyclicity Criterion is that it is actually equivalent to an operator being weakly mixing. 
In \cite{BP99}, B\`es and Peris obtained this result using the notion of hereditarily hypercyclic operators.

\begin{thm}[{\cite[Theorem 2.3]{BP99}}]\label{BP}
Let $T\colon X\to X$ be an operator. 
Then the following assertions are equivalent:
    \begin{enumerate}
        \item $T$ satisfies the Hypercyclicity Criterion;
        \item $T$ is weakly mixing;
        \item $T$ is hereditarily hypercyclic.
    \end{enumerate}
\end{thm}

Unlike the case for compact dynamical systems, a tool for the study of linear dynamical systems is provided by the neighborhoods of $0$. 
We recall the Blow up/Collapse condition for weak mixing.

\begin{thm}[{\cite[Theorem 2.47]{G11}}] \label{blowup-T}
Let $T\colon X\to X$ be an operator. Then $T$ is weakly mixing if and only if for any non-empty open sets $U, V\subset X$ and any neighborhood $W$ of $0$,
     \[
         N(U, W)\cap N(W, V)\neq\emptyset. 
     \]
\end{thm}

Let $X= c_0(\mathbb{N}_0)$ or $\ell^p(\mathbb{N}_0)$, $1\leq p<\infty$, and let $w=(w_i)_{i\in\mathbb{N}_0}$ be a bounded sequence of positive real numbers. 
Define a map $B_w\colon X\to X$ by 
    \[
        (x_0, x_1, x_2, \dotsc )\mapsto (w_1x_1, w_2x_2, \dotsc).
    \]
The map $B_w$ is clearly an operator on $X$, and is called a \emph{unilateral weighted backward shift}. 
Similarly, one can define \emph{bilateral weighted backward shift} on $c_0(\mathbb{Z})$ and $\ell^p(\mathbb{Z})$, $1\leq p<\infty$.
In \cite{BMPP19}, B\`es et al.\@ characterized the $\mathcal{F}$-transitivity of such operators, and one of their results is restated here. 

\begin{prop}[{\cite[Proposition 3.1]{BMPP19}}]\label{Bes-char}
Let $\mathcal{F}$ be a family and $B_w$ be a bilateral weighted backward shift on $X=c_0(\mathbb{Z})$ or $\ell^p(\mathbb{Z})$, $1\leq p<\infty$. 
Then the following assertions are equivalent: 
    \begin{enumerate}
        \item $B_w$ is $\mathcal{F}$-transitive; 
        \item $B_w$ is $\mathcal{F}$-mixing (i.e.\ $\mathcal{F}$-transitive and weakly mixing); 
        \item for any $M>0$ and $j\in\mathbb{N}$, we have $A_{M,j}\cap \overline{A}_{M,j}\in\mathcal{F}$, where 
        \[
            A_{M,j}=\{n\in\mathbb{N}\colon w_{j+1}w_{j+2}\dotsb w_{j+n}>M\}, \overline{A}_{M,j}=\{n\in\mathbb{N} \colon  w_{j}w_{j-1}\dotsb w_{j-n+1}<\tfrac{1}{M}\}.
        \]
    \end{enumerate}
\end{prop}

We denote by $\omega$ the Fr\'echet space of all sequences in $\mathbb{K}$, where $\mathbb{K}$ denotes either  $\mathbb{R}$ or $\mathbb{C}$, that is, 
\[
    \omega \colon=\mathbb{K}^{\mathbb{N}_0}=\{(x_n)_n\colon x_n\in \mathbb{K}, n\in \mathbb{N}_0\},
\]
with the topology given by the seminorms: $p_n(x)=\operatorname{sup}_{i\leq n}|x_i|, n\in\mathbb{N}_0$. 
The \emph{canonical basis} of $\omega$ is denoted by $e_n=(\delta_{i,n})_{i\in\mathbb{N}_0}$, $n\in\mathbb{N}_0$. 

Let $X$ be an infinite-dimensional separable Fr\'echet space. 
If $X$ can be embedded continuously (by a linear map, $\theta$ say) into $\omega$, then $X$ is called a \emph{unilateral Fr\'echet sequence space}. 
Roughly speaking, the convergence in a Fr\'echet sequence space implies the convergence of each coefficient in $\mathbb{K}$. 
Since for each $n\in\mathbb{N}_0$, $E_n=\theta^{-1}e_n\in X$ (if exists) must be unique, we denote it still by $e_n$. 

Let $X$ be a unilateral Fr\'echet sequence space in which $(e_n)_{n\in\mathbb{N}_0}$ is a basis (i.e.\@ $\sum_{i=0}^{\infty}x_ie_i$ converges for each $x\in X$). Define a linear map $B\colon X\to \omega$ by
    \[
        (x_0, x_1, x_2, \dotsc )\mapsto (x_1, x_2, \dotsc),
    \]
where $x=\sum_{i=0}^{\infty}x_ie_i\in X.$
If $B$ maps $X$ into $X$, then it must be continuous and we say that $B$ is the \emph{unilateral backward shift} on $X$. 
Define a map $F\colon \{x\in X\colon \text{there is }l\in\mathbb{N}_0\text{ such that }x_i=0\text{ for each }i\geq l\}\to X$ by 
    \[
        (x_0, x_1, x_2, \dotsc,x_l,0\dotsc )\mapsto (0,x_0, x_1, x_2, \dotsc,x_l,0\dotsc  ).
    \]
The map $F$ is called the \emph{unilateral forward shift for finite sequences} of $X$. 
\emph{Bilateral backward (and forward) shifts} on \emph{bilateral Fr\'echet sequence spaces} are defined similarly. 
See \cite{BMPP16}, \cite{BMPP19}, \cite{G11} and \cite{S95} for more details of weighted backward shifts.

\section{Weak disjointness in linear dynamics}

In this section we introduce the notion of weak disjointness of linear dynamical systems. 
We will give a criterion for $\mathcal{F}$-mixing, and a sufficient condition for mild mixing. 
Also, we illustrate that a weakly mixing operator may not be weakly disjoint from its iterates. 

\begin{defn}
Let $T\colon X\to X$ and $S\colon Y\to Y$ be two operators. 
We say that $T$ and $S$ are \emph{weakly disjoint} provided that their direct sum $T\oplus S$ is hypercyclic. 
\end{defn}

Weak disjointness of finitely many operators is defined similarly. 
By Remark \ref{Birkhoff}, two operators $T\colon X\to X$ and $S\colon Y\to Y$ are weakly disjoint if and only if for any non-empty open subsets $U_1,U_2$ of $X$ and non-empty open subsets $V_1,V_2$ of $Y$, $N_T(U_1,U_2)\cap N_S(V_1,V_2)\neq\emptyset$. 
We have the following observation, which is in fact true for topological dynamical systems. 
Recall that for a given family $\mathcal{F}$, an operator $T\colon X\to X$ is $\mathcal{F}$-transitive if $N(U,V)\in\mathcal{F}$ for any non-empty open subsets $U$ and $V$ of $X$.
    
\begin{lem}\label{F-disjoint}
Suppose that $\mathcal{F}$ is a shift$_+$-invariant family. 
Let $T\colon X\to X$ and $S\colon Y\to Y$ be two operators. 
If $S$ is hypercyclic, then $T\oplus S$ is $\mathcal{F}$-transitive if and only if for any non-empty open subsets $U_1,U_2$ of $X$ and  non-empty open subset  $V$ of $Y$, we have
\[
    N_T(U_1,U_2)\cap N_S(V,V)\in\mathcal{F}.
\]
\end{lem}

\begin{proof}
The condition is obviously necessary.

Conversely, by the hypothesis, both $T$ and $S$ are hypercyclic. Fix non-empty open subsets $U_1,U_2$ of $X$ and  non-empty open subsets $V_1,V_2$ of $Y$. By the hypercyclicity of $S$, there must be $m\in\mathbb{N}_0$ such that $V_1\cap S^{-m}V_2\neq\emptyset$. By the hypercyclicity of $T$, $T^{-m}U_2\neq\emptyset$. Let $V=V_1\cap S^{-m}V_2$. We then have $N_S(V, V)+m\subset N_S(V_1, V_2)$ and $N_T(U_1, T^{-m}U_2)+m\subset N_T(U_1, U_2)$.
Finally
\[
    (N_T(U_1, T^{-m}U_2)\cap N_S(V, V))+m\subset N_T(U_1,U_2)\cap N_S(V_1,V_2).
\]
We are done since $N_T(U_1, T^{-m}U_2)\cap N_S(V, V)\in\mathcal{F}$ and $\mathcal{F}$ is shift$_+$-invariant. 
\end{proof}

Since the family of non-empty subsets of $\mathbb{N}_0$ is shift$_+$-invariant, applying Lemma \ref{F-disjoint} we have the following version for weak disjointness. 
    
\begin{lem}\label{center}
Let $T\colon X\to X$ and $S\colon Y\to Y$ be two operators. 
If $S$ is hypercyclic, then $T$ and $S$ are weakly disjoint if and only if for any non-empty open subsets $U_1,U_2$ of $X$ and  non-empty open subset  $V$ of $Y$, we have $$N_T(U_1,U_2)\cap N_S(V,V)\neq\emptyset. $$
\end{lem}

Given a hypercyclic operator $T$ and $k\in\mathbb{N}$, we consider the weak disjointness of $T$ and $T^k$. 
We will see in the following examples that $T$ may not be weakly disjoint from $T^k$, even if $T$ is weakly mixing. 
We note that related constructions appear in \cite[Theorem 1.3]{GR09} and \cite[Example 4.5]{BP07}, although they are used there to address different questions. 
We give a self-contained proof below, formulated in a way that makes explicit the properties required for our subsequent analysis. 

\begin{exam}\label{example1}
Let $X=c_0(\mathbb{N}_0)$ or $\ell^p(\mathbb{N}_0)$, $1\leq p<\infty$. 
There exists a unilateral weighted backward shift $B_w$ on $X$ such that $B_w$ is weakly mixing but $B_w\oplus B^2_w$ is not hypercyclic. 
\end{exam}

\begin{proof}
For any $n\in\mathbb{N}$, we define two finite blocks $\mathcal{A}_n$ and $\mathcal{B}_n$ by 
    \[
        \mathcal{A}_n=[\underbrace{2,2,2,\dotsc, 2}_{n\text{ terms}}, \tfrac{1}{2^n}],\ \mathcal{B}_n=[\underbrace{1,1,1,\dotsc, 1}_{n\text{ terms}}].
    \]
Define a sequence $(b_n)_{n\in\mathbb{N}}$ of positive integers by 
\[
    b_1=2,\ b_{n+1}=2b_n+n+2,\quad n\in\mathbb{N}. 
\]
Let $B_w$ be the unilateral weighted backward shift on $X$ defined by the weight $w$, where $w=(w_i)_{i\in\mathbb{N}_0}$ is given by 
\[
        w=(\mathcal{B}_1,\underbrace{\mathcal{A}_1, \mathcal{B}_{b_1}}, \underbrace{\mathcal{A}_2,\mathcal{B}_{b_2}},\dotsc,\underbrace{\mathcal{A}_n,\mathcal{B}_{b_n}},\dotsc). 
\]
For a finite block $u$, denote the length of $u$ by $|u|$. 
For each $n\in\mathbb{N}$, one has 
\[ 
\sum_{i=1}^n (|\mathcal{A}_i|+|\mathcal{B}_{b_i}|)=2b_n.
\]
For each $n\in\mathbb{N}$, it is easy to see that for each $k\in [b_n-1,2b_n]\cap\mathbb{N}_0$, 
$w_1w_2\dotsb w_k=1$,
and for each $k\in [2b_n+1,2b_n+n+1]\cap\mathbb{N}_0$,
$w_1w_2\dotsb w_k=2^{k-2b_n}$.
Then $\limsup_{k\to\infty}w_1w_2\dotsb w_{k}=\infty$. 
By \cite[Example 4.9]{G11}, $B_w$ is weakly mixing. 

For each $k\in\mathbb{N}$ with $w_1w_2\dotsb w_k\geq 2$, we have either $k=1$ or $k\in [2b_n+1,2b_n+n+1]$ for some $n\in\mathbb{N}$. 
If $k=1$ then $w_1w_2\dotsb w_{2k}=w_1w_2=1$. 
If $k\in [2b_n+1,2b_n+n+1]$ for some $n\in\mathbb{N}$, 
since $2k\in [4b_n+2, 4b_n+2n+2]\subset [b_{n+1}-1,2b_{n+1}]$, 
we have $w_1w_2\dotsb w_{2k}=1$. 
Consider the non-empty open set 
\[
    W=\bigl\{x\in X\colon \|x-e_0\|<\tfrac{1}{3}\bigr\}.
\]
Assume that there is $m\in N_{B_w\oplus B_w^2}(W\oplus W, W\oplus W)\setminus\{0\}$. 
Then there is a pair $(x,y)\in W\oplus W$ with 
$(B_w\oplus B_w^2)^m (x,y)\in W\oplus W$.
Since $x,y\in W$ and $m\neq 0$, we have $|x_m|<\frac{1}{3}$ and $|y_{2m}|<\frac{1}{3}$.
Since $B_w^m x\in W$, we have $|w_1w_2\dotsb w_m x_m-1|<\frac{1}{3}$, and then $w_1w_2\dotsb w_m=\tfrac{w_1w_2\dotsb w_m|x_m|}{|x_m|}>2$. 
Then $w_1w_2\dotsb w_{2m}=1$ by the above analysis. 
Since $B_w^{2m} y\in W$, we have $|w_1w_2\dotsb w_{2m} y_{2m}-1|<\frac{1}{3}$, and then again $w_1w_2\dotsb w_{2m}>2$, which is a contradiction.
Then $N_{B_w\oplus B_w^2}(W\oplus W, W\oplus W)=\{0\}$, which implies that $B_w\oplus B^2_w$ is not hypercyclic. 
\end{proof}

\begin{exam}\label{B-B^k}
Let $X=c_0(\mathbb{N}_0)$ or $\ell^p(\mathbb{N}_0)$, $1\leq p<\infty$. 
There exists a unilateral weighted backward shift $B_w$ on $X$ such that $B_w\oplus B^2_w$ is weakly mixing but $B_w\oplus B^3_w$ is not hypercyclic. 
\end{exam}

\begin{proof}
Assume that $X=c_0(\mathbb{N}_0)$, the case $X=\ell^p(\mathbb{N}_0)$, $1\leq p<\infty$ is similar. 
For any $n, m\in\mathbb{N}$, we define three finite blocks $\mathcal{A}_{n},\mathcal{B}_{n,m}$ and $\mathcal{C}_n$ by 
    \[
        \mathcal{A}_{n}=[\underbrace{2,2,2,\dotsc, 2}_{n\text{ terms}}],\ \mathcal{B}_{n,m}=[\underbrace{1,1,1,\dotsc, 1}_{n\text{ terms}},\tfrac{1}{2^m}],\ \mathcal{C}_n=[\underbrace{1,1,1,\dotsc, 1}_{n\text{ terms}}]. 
    \]
Let $s_0=1,s_1=3$ and $s_{n+1}=2(3(s_n-1)+n+1)+n$ for $n\in\mathbb{N}$. 
Define two sequences $(b_n)_{n\in\mathbb{N}}$ and $(c_n)_{n\in\mathbb{N}}$ of positive integers by 
    \[
        b_n=s_n-1-[3(s_{n-1}-1)+n],\ c_n=2(s_n-1)-1, \quad n\in\mathbb{N}. 
    \]
Let $B_w$ be the unilateral weighted backward shift on $X$ defined by the weight $w$, where $w=(w_i)_{i\in\mathbb{N}_0}$ is given by 
    \[
        w=(\mathcal{C}_1,\underbrace{\mathcal{A}_1,\mathcal{B}_{b_1,1}, \mathcal{C}_{c_1}},\dotsc,\underbrace{\mathcal{A}_n,\mathcal{B}_{b_n,n}, \mathcal{C}_{c_n}},\dotsc). 
    \]
Fix $n\in\mathbb{N}$. 
One has $w_{s_n}=2^{-n}$. 
If $k\in[s_n,3(s_n-1)]\cap\mathbb{N}_0$ then $w_1w_2\dotsb w_k=1$. 
If $k\in[3(s_n-1)+1,3(s_n-1)+n]\cap\mathbb{N}_0$ then $w_1w_2\dotsb w_k=2^{k-3(s_n-1)}$. 
If $k\in[3(s_n-1)+n+1,s_{n+1}-1]\cap\mathbb{N}_0$ then $w_1w_2\dotsb w_k=2^{n+1}$. 

For each finite sequence $x=(x_0,x_1,x_2,\dotsc,x_l,0,\dotsc)\in X$, let 
\[
    S_wx=(0,\tfrac{x_0}{w_1},\tfrac{x_1}{w_2},\tfrac{x_2}{w_3},\dotsc,\tfrac{x_l}{w_{l+1}},0,\dotsc). 
\]
Denote $T=B_w\oplus B_w^2$, $S=S_w\oplus S^2_w$ and $Y=X\oplus X$. 
Let $X_0$ be the set of finite sequences in $X$. 
Then $Y_0=X_0\oplus X_0$ is a dense subset of $Y$. 
Put $a_n=3(s_{n-1}-1)+n$ for each $n\in\mathbb{N}$. 
We shall show that $T$ satisfies the Hypercyclicity Criterion (Theorem \ref{HC}) with $S$, $Y_0$ and $(a_n)_{n\in\mathbb{N}}$. 

To see this, fix $(u,v)\in Y_0$ and $\epsilon>0$. 
Write $u=(u_0,u_1,\dotsc, u_l,0,0,\dotsc)$ and $v=(v_0,v_1,\dotsc, v_l,0,0,\dotsc)$. 
Let $M_1=\operatorname{max}(\{|u_i|: i\in\mathbb{N}_0\}\cup \{|v_i|: i\in\mathbb{N}_0\})+1, M_2=\operatorname{max}\{w_1w_2\cdots w_i: i\in\mathbb{N}_0, i\leq l\}$. 
Pick $N_1\in\mathbb{N}$ such that $n>N_1$ implies $s_n-2n-7>l$. 
For any $n>N_1$ we have $l+2a_n-a_n<s_{n}-1-(3(s_{n-1}-1)+n)$ and then $\{a_n,a_n+1,\dotsc, 2a_n, \dotsc,l+2a_n\}\subset [3(s_{n-1}-1)+n,s_{n}-1]$. 
By the above analysis, we have $w_{1}\dotsb w_{a_n}=w_{1}\dotsb w_{a_n+1}=\dotsc =w_{1}\dotsb w_{l+2a_n}=2^n$. 
Pick $N_2\in\mathbb{N}$ such that $n>N_2$ implies $\tfrac{1}{2^n}<\tfrac{\epsilon}{M_1M_2}$. Put $N=\operatorname{max}\{N_1,N_2\}$. 
For any $n>N$, one has 
\begin{align*}
    \|S_w^{a_n}u\|&=\operatorname{sup}_{j\leq l} \frac{|u_{j}|}{w_{j+1}\cdots w_{j+a_n}}
        =\operatorname{sup}_{j\leq l} \frac{w_1w_2\dotsb w_j|u_{j}|}{w_1w_2\dotsb w_jw_{j+1}\cdots w_{j+a_n}}\\
        &=\operatorname{sup}_{j\leq l} \frac{w_1w_2\dotsb w_j|u_{j}|}{2^n}
        <\operatorname{sup}_{j\leq l} \frac{w_1w_2\dotsb w_j|u_{j}|}{M_1M_2}\epsilon\\
        &<\epsilon, 
\end{align*}
and 
\begin{align*}
    \|S_w^{2a_n}v\|&=\operatorname{sup}_{j\leq l} \frac{|v_{j}|}{w_{j+1}\cdots w_{j+2a_n}}
        =\operatorname{sup}_{j\leq l} \frac{w_1w_2\dotsb w_j|v_{j}|}{w_1w_2\dotsb w_jw_{j+1}\cdots w_{j+2a_n}}\\
        &=\operatorname{sup}_{j\leq l} \frac{w_1w_2\dotsb w_j|v_{j}|}{2^n}
        <\operatorname{sup}_{j\leq l} \frac{w_1w_2\dotsb w_j|v_{j}|}{M_1M_2}\epsilon\\
        &<\epsilon. 
\end{align*}
This gives that $\displaystyle\lim_{n\to\infty}S^{a_n}(u,v)=0$, while it is clear that $TS(u,v)=(u,v)$ and $\displaystyle\lim_{n\to\infty}T^{a_n}y_0=0$ for each $y_0\in Y_0$. 
So $T$ satisfies the Hypercyclicity Criterion (Theorem \ref{HC}) and then it is weakly mixing. 

But whenever there is $n\in\mathbb{N}$ with $w_1w_2\cdots w_n\geq 2$, we have either $n\in \{1,2\}$ or $n\in [3(s_k-1)+1,s_{k+1}-1]$ for some $k\in \mathbb{N}$. 
If $n\in \{1,2\}$ then $w_1w_2\cdots w_{3n}=1$. 
If $n\in [3(s_k-1)+1,s_{k+1}-1]$ for some $k\in \mathbb{N}$, since $3n\in [9(s_k-1)+3,3(s_{k+1}-1)]\subset [s_{k+1},3(s_{k+1}-1)]$, we have $w_1w_2\cdots w_{3n}=1$. 
Consider the non-empty open set 
\[
    W=\bigl\{x\in X\colon \|x-e_0\|<\tfrac{1}{3}\bigr\}.
\]
As in Example \ref{example1}, $N_{B_w\oplus B_w^3}(W\oplus W, W\oplus W)=\{0\}$ and then $B_w\oplus B_w^3$ is not hypercyclic. 
\end{proof}

\begin{rem}
By a similar method, for any $k\in\mathbb{N}$ and $X=c_0(\mathbb{N}_0)$ 
or $\ell^p(\mathbb{N}_0)$, $1\leq p<\infty$, one can construct a unilateral weighted backward shift $B_w$ on $X$ such that $B_w\oplus B_w^2\oplus\dotsb \oplus B^k_w$ is weakly mixing but $B_w\oplus B_w^{k+1}$ is not hypercyclic. 
\end{rem}
   
\subsection{The \texorpdfstring{$\mathcal{F}$}
{mathcal{F}}-Mixing Criterion}

It is observed that an operator is weakly mixing if and only if it is weakly disjoint from itself. 
We investigate how a weakly mixing operator can be weakly disjoint from itself.
In terms of Furstenberg families, we study the $\mathcal{F}$-transitivity of the direct sum of an operator and itself when a family $\mathcal{F}$ is given. 

\begin{defn}
    Let $\mathcal{F}$ be a family. 
    An operator $T\colon X\to X$ is said to be \emph{$\mathcal{F}$-mixing} provided that $T\oplus T$ is $\mathcal{F}$-transitive. 
\end{defn}

Recall that the Blow up/Collapse condition characterizes the weak mixing of operators (see Theorem \ref{blowup-T}). 
We formulate here the Blow up/Collapse condition for $\mathcal{F}$-mixing. 
    
\begin{lem}\label{F-mixing-blow-up}
    Let $\mathcal{F}$ be a family. 
    Then an operator $T\colon X\to X$ is $\mathcal{F}$-mixing if and only if for any non-empty open subsets $U,V$ of $X$ and neighborhood $W$ of \ $0$, we have 
    \[
        N(U,W)\cap N(W,V)\in\mathcal{F}. 
     \]
\end{lem}

\begin{proof}
The condition is obviously necessary.

Conversely, let $U,V$ be a pair of non-empty open subsets of $X$. 
Pick non-empty open subsets $U'\subset U,V'\subset V$ and a neighborhood $W$ of $0$ such that $U'+W\subset U, V'+W\subset V$. 
The observation $N(U',W)\cap N(W,V')\subset N(U,V)$ gives that $N(U,V)\in\mathcal{F}$ and then $T$ is $\mathcal{F}$-transitive. 
Since $\emptyset\notin\mathcal{F}$, by Theorem \ref{blowup-T}, $T$ is weakly mixing. 
So $T$ is $\mathcal{F}$-mixing. 
\end{proof}

Inspired by the Hypercyclicity Criterion (Theorem \ref{HC}), we formulate a criterion for $\mathcal{F}$-mixing in terms of limits, which is especially applicable to backward shifts (see Section 5). 
Before this, we need a notion, which requires the specific shape of the indices of convergent subsequences. 
    
Suppose that $X$ is a metric space and $\mathcal{F}$ is a family. 
Let $(x_n)_{n\in\mathbb{N}_0}$ be a sequence in $X$. 
We say that a point $x\in X$ is an \emph{$\mathcal{F}$-limit} of $(x_n)_{n\in\mathbb{N}_0}$, denoted by $\mathcal{F}$-$\displaystyle\lim_{n\to\infty}x_n=x$, if for any neighborhood $U$ of $x$ we have 
\[
\{n\in\mathbb{N}_0: x_n\in U\}\in\mathcal{F}. 
\]

In \cite{BMPP19}, B\`es et al.\@ gave a criterion \cite[Theorem 2.4]{BMPP19} for an operator to be $\mathcal{F}$-mixing, requiring that $\widetilde{\mathcal{F}}$ is a filter. 
The following criterion has a strong requirement on the limits but no requirement for the family. 

\begin{thm}[$\mathcal{F}$-Mixing Criterion]\label{F-mixing1}
Let $T\colon X\to X$ be an operator and $\mathcal{F}$ be a family.
If there exist dense subsets $X_0$, $Y_0$ of $X$ and a sequence $(S_n)_{n\in\mathbb{N}_0}$ of maps of $Y_0$ into $X$ such that for any $x_0\in X_0$, $y_0\in Y_0$, 
    \[
        \mathcal{F}\text{-}\displaystyle\lim_{n\to\infty}(T^{n}x_0, S_{n}y_0, T^{n}S_{n}y_0)=(0, 0, y_0), 
    \]
then $T$ is $\mathcal{F}$-mixing. 
\end{thm}

\begin{proof}
Let $U, V$ be non-empty open subsets of $X$ and let $W\subset X$ be a neighborhood of $0$. 
Pick $x_0\in X_0\cap U$ and $y_0\in Y_0\cap V$. 
Then 
\[
A=\{n\in\mathbb{N}_0: T^nx_0\in W, S_ny_0\in W, T^nS_ny_0\in V\}\in\mathcal{F}. 
\]
Since $A\subset N(U, W)\cap N(W, V)$, we have $N(U,W)\cap N(W,V)\in\mathcal{F}$. Then $T$ is $\mathcal{F}$-mixing, by Lemma \ref{F-mixing-blow-up}. 
\end{proof}

\begin{rem}
It is clear that an operator satisfying the Hypercyclicity Criterion (Theorem \ref{HC}) must satisfy the $\mathcal{F}_{\mathrm{inf}}$-Mixing Criterion. 
Once an operator satisfies the $\mathcal{F}_{\mathrm{inf}}$-Mixing Criterion, it is weakly mixing and then by Theorem \ref{BP} it satisfies the Hypercyclicity Criterion (Theorem \ref{HC}). 
The $\mathcal{F}_{\mathrm{cf}}$-Mixing Criterion coincides with Kitai's Criterion (\cite[Theorem 3.4]{G11}). 
Generally, the $\mathcal{F}$-Mixing Criterion does not characterize the $\mathcal{F}$-mixing. 
For example, Grivaux has constructed a strongly mixing operator which does not satisfy Kitai's Criterion (see \cite[Theorem 2.5]{G05}). 
\end{rem}

\subsection{Mild mixing}
It is natural to consider the property of operators which are weakly disjoint from a given class of operators.
In 2024, Cardeccia characterized in \cite{C24} the operators which are weakly disjoint from all weakly mixing operators, and the operators which are weakly disjoint from all topologically ergodic operators. We restate these results, in the language of weak disjointness. 

\begin{thm}[{\cite[Theorem 5.3]{C24}}]\label{wk-mixing-dsj}
Let $T\colon X\to X$ be an operator. 
The following assertions are equivalent: 
    \begin{enumerate}
        \item $T$ is topologically ergodic; 
         \item $T$ is weakly disjoint from all weakly mixing operators; 
        \item $T\oplus S$ is weakly mixing for each weakly mixing operator $S\colon Y\to Y$. 
    \end{enumerate}
\end{thm}

\begin{thm}[{\cite[Theorem 5.4]{C24}}]\label{ps-dsj}
An operator $T\colon X\to X$ is $\mathcal{F}_{\mathrm{ps}}$-transitive if and only if $T$ is weakly disjoint from all topologically ergodic operators.
\end{thm}

It is clear that the operators which are weakly disjoint from all strongly mixing operators are exactly the hypercyclic ones.
We wonder if there is a non-mixing operator which is weakly disjoint from all hypercyclic operators.
We introduce the concept of mild mixing in the framework of linear dynamical systems, and show that the mild mixing of operators is a property which lies strictly between strong mixing and weak mixing. 

\begin{defn}
An operator $T\colon X\to X$ is said to be \emph{mildly mixing} provided that it is weakly disjoint from all hypercyclic operators. 
\end{defn}

\begin{prop}
Let $T\colon X\to X$ and $S\colon Y\to Y$ be two operators. Then $T\oplus S$ is mildly mixing if and only if both $T$ and $S$ are mildly mixing. 
\end{prop}

\begin{proof}
The condition is obviously necessary.

Conversely, fix a hypercyclic operator $R\colon Z\to Z$. 
Since $S$ is mildly mixing, we have that $S\oplus R$ is hypercyclic. 
Since $T$ is mildly mixing, 
we have that $T\oplus (S\oplus R)$ is hypercyclic. This gives that $(T\oplus S)\oplus R$ is hypercyclic and then $T\oplus S$ is mildly mixing. 
\end{proof}

Next we deduce a sufficient condition for an operator to be mildly mixing. 

\begin{prop}\label{mild-ip}
Each $\Delta^*_{\mathrm{ip}}$-transitive operator is mildly mixing. 
\end{prop}

\begin{proof}
Suppose that $T\colon X\to X$ is a $\Delta^*_{\mathrm{ip}}$-transitive operator. 
Fix a hypercyclic operator $S\colon Y\to Y$. 
Let $U_1, U_2$ be non-empty open subsets of $X$ and let $V$ be a non-empty open subset of $Y$. 
Let $x\in V$ be a hypercyclic vector for $S$. 
Since each hypercyclic vector is recurrent, by Theorem~\ref{NxW-ip}, we have $N_S(x, V)\in\mathcal{F}_{\mathrm{ip}}$ and $N_S(V,V)=(N_S(x,V)-N_S(x, V))\cap \mathbb{N}_0$ (i.e.\@ $N_S(V,V)\in \Delta_{\mathrm{ip}}$). 
Since $N_T(U_1, U_2)\in\Delta^*_{\mathrm{ip}}$, we have that $N_T(U_1, U_2)\cap N_S(V,V)\neq\emptyset$.
By Lemma \ref{center}, $T$ is weakly disjoint from $S$ and then $T$ is mildly mixing. 
\end{proof}

Recall from Theorem \ref{cpt-mm-ip} we know that the mild mixing of compact dynamical systems is equivalent to the $\Delta^*_{\mathrm{ip}}$-transitivity. 
In Proposition \ref{mild-ip} we have shown that each $\Delta^*_{\mathrm{ip}}$-transitive operator is mildly mixing, so we ask: 

\begin{Ques}
Is each mildly mixing operator $\Delta^*_{\mathrm{ip}}$-transitive? 
\end{Ques}

Now we deduce that for operators, strong mixing is not equivalent to mild mixing. 
Let $\Delta=\{A\subset\mathbb{N}_0: \text{there is }F\in\mathcal{F}_{\mathrm{inf}}\text{ with }(F-F)\cap\mathbb{N}_0\subset A\}$. 
Since $\mathcal{F}_{\mathrm{ip}}\subset\mathcal{F}_{\mathrm{inf}}$, by Proposition \ref{mild-ip}, each $\Delta^*$-transitive operator is mildly mixing. 
B\`es et al.\@ have constructed a $\Delta^*$-transitive weighted backward shift on $c_0(\mathbb{N}_0)$ or $\ell^p(\mathbb{N}_0), 1\leq p<\infty$ which is not strongly mixing in \cite[Proposition 5.5]{BMPP19}. So we have the following. 

\begin{prop}\label{mm-sm}
There exists a mildly mixing operator which is not strongly mixing. 
\end{prop}

\begin{rem}\label{TT^*}
For an operator $T$ on $\ell^2(\mathbb{Z})$, it is natural to consider the weak disjointness of $T$ and its adjoint $T^*$. 
But it can be easily checked that this will never happen, see e.g.\@ \cite[Remark 4.17]{G11}. 
\end{rem} 

Salas constructed in \cite[Corollary 2.6]{S95} two weakly mixing operators $T,S$ on $\ell^2(\mathbb{Z})$ such that $T$ and $S$ are not weakly disjoint.
He also constructed in \cite[Corollary 2.3]{S95} an operator $T$ on $\ell^2(\mathbb{Z})$ such that both $T$ and $T^*$ are weakly mixing, while $T$ and $T^*$ can never be weakly disjoint by Remark \ref{TT^*}. 
It is illustrated in Example \ref{B-B^k} that there is a weakly mixing operator $T$ on $\ell^1(\mathbb{N}_0)$ such that $T$ and $T^2$ are not weakly disjoint. 
The above examples show that mild mixing is not equivalent to weak mixing for operators. 

\begin{prop}\label{mm-wm}
There exist some weakly mixing operators which are not mildly mixing. 
\end{prop}

It is observed that each strongly mixing operator is mildly mixing and each mildly mixing operator is weakly mixing. 
By Propositions \ref{mm-sm} and \ref{mm-wm}, the mild mixing of operators lies strictly between strong mixing and weak mixing. 
In 1993, Ansari showed in \cite{A95} that an operator is hypercyclic if and only if each of its iterates is hypercyclic. 

\begin{thm}[{\cite[Theorem 1]{A95}}]\label{Ansari}
Let $T\colon X\to X$ be an operator and $p\in\mathbb{N}$. 
Then $HC(T)=HC(T^p)$. 
\end{thm}

By Ansari's result, we have the following two more sufficient conditions for mild mixing. 
Recall that an operator is said to be \emph{chaotic in the sense of Devaney} (or \emph{chaotic} for short) if it is hypercyclic and has a dense set of periodic points. 
In 2007, Badea and Grivaux constructed in \cite[Corollary 4.7]{BG07} a chaotic operator which is not strongly mixing, which shows that chaos does not imply strong mixing. 
But we have that chaos does imply mild mixing.
We also obtain that an operator is mildly mixing if some iterate of it is mildly mixing. 
 
\begin{prop}
Let $T\colon X\to X$ be an operator. 
\begin{enumerate}
    \item If $T$ is chaotic, then it is mildly mixing; 
    \item If $T^p$ is mildly mixing for some $p\in\mathbb{N}$, then $T$ is mildly mixing. 
\end{enumerate}
\end{prop}

\begin{proof}
(1) Let $S\colon Y\to Y$ be a hypercyclic operator. 
Fix a non-empty open subset $V$ of $X$ and two non-empty open subsets $U_1, U_2$ of $Y$. 
Since $T$ is chaotic, there is $v\in V$ and $p\in\mathbb{N}$ such that $T^pv=v$. 
 Since $S^p$ is hypercyclic by Theorem \ref{Ansari}, there exists $n\in\mathbb{N}$ such that $S^{pn}U_1\cap U_2\neq\emptyset$. 
Since $T^{pn}v=v$, we have $pn\in N_T(V,V)\cap N_S(U_1,U_2)$. 
By Lemma \ref{center}, $T$ is weakly disjoint from $S$. 
Thus $T$ is mildly mixing. 

(2) Let $S\colon Y\to Y$ be a hypercyclic operator. 
By Theorem \ref{Ansari}, $S^p$ is also hypercyclic. 
Since $T^p$ is mildly mixing, 
we have that $T^p$ is weakly disjoint from $S^p$.
Then $T$ is weakly disjoint from $S$.
This implies that $T$ is mildly mixing.
\end{proof}

It is natural to consider the following question.
\begin{Ques}
Let $T\colon X\to X$ be a mildly mixing operator.  
Is $T^p$ mildly mixing for all $p\in\mathbb{N}$?
\end{Ques}

\section{The Weiss-Akin-Glasner Theorem for operators and some applications}

In this section, by showing an analogue of the Weiss-Akin-Glasner Lemma, we formulate the counterpart of the Weiss-Akin-Glasner Theorem within the framework of linear dynamical systems, and provide some applications. 

\subsection{The Weiss-Akin-Glasner Theorem for operators}
We have the following  characterization of subsets of non-negative integers which can be realized by the return time sets of $\mathcal{F}$-mixing operators,
which is an analogue of the Weiss-Akin-Glasner Lemma (Lemma~\ref{lem:WAG-revised}) from topological dynamics.

\begin{lem}\label{L-WAG}
Let $\mathcal{F}$ be a full family and $F\subset \mathbb{N}_0$.
Then
$F\in \widetilde{\mathcal{F}}$ if and only if 
there is an invertible $\mathcal{F}$-mixing weighted backward shift $B_w$ on $X$, and a non-empty open subset $W$ of $X$ such that 
\[
    N_{B_w}(W,W)=F\cup\{0\}, 
\]
where $X=c_0(\mathbb{Z})$ or $ \ell^p(\mathbb{Z})$, $1\leq p<\infty$.
\end{lem}

\begin{proof}
If such an operator $B_w\colon X\to X$ exists, then by Proposition~\ref{return-time-F-mixing}, 
$\mathcal{N}_{B_w}$ is a shift-invariant thick subfamily of $\mathcal{F}$. 
Then $F\cup\{0\} \in \mathcal{N}_{B_w}\subset \widetilde{\mathcal{F}}$ since $\widetilde{\mathcal{F}}$ is the maximal shift-invariant thick subfamily of $\mathcal{F}$. 
Since $\widetilde{\mathcal{F}}$ is shift-invariant, it is full. 
Then $(F\cup\{0\})\cap \mathbb{N}\in\widetilde{\mathcal{F}}$. 
It is observed that $(F\cup\{0\})\cap \mathbb{N}\subset F$. 
So $F\in\widetilde{\mathcal{F}}$.

Conversely, suppose that $F\in\widetilde{\mathcal{F}}$. 
Since $\widetilde{\mathcal{F}}$ is hereditary upward, $F\cup\{0\}\in \widetilde{\mathcal{F}}$. 
As the desired return time set is $F\cup\{0\}$, replacing $F$ by $F\cup\{0\}$ if necessary, we may assume that $0\in F$. 
We divide the construction of the weighted backward shift into two cases. 

\textbf{Case A: }$F\in \mathcal{F}_{\mathrm{cf}}$. 
Write $F=I_1\cup I_2\cup\dotsb \cup I_t\cup ([L,\infty)\cap\mathbb N_0)$ with $t,L\geq 0$ such that $\{I_k\}_{k=1}^t$ are disjoint maximal finite integer intervals. 
Write $\mathbb{N}_0\setminus F=J_1\cup J_2\cup\dotsb \cup J_t$ such that $\{J_k\}_{k=1}^t$ are disjoint maximal finite integer intervals. 
Without loss of generality we assume that $\min I_1<\min I_2<\dotsb<\min I_t$ and $\min J_1<\min J_2<\dotsb<\min J_t$. 
If $F=\mathbb{N}_0$, then there is no such $I_k$ and $J_k$.
For each $k\in\{1,2,\dotsc, t\}$, let
\begin{align*} 
\mathcal{I}_k=
    \begin{cases}
    [\underbrace{2,2,\dotsc,2}_{s\text{ terms}},\underbrace{\tfrac{1}{2},\tfrac{1}{2},\dotsc,\tfrac{1}{2}}_{s-1\text{ terms}},1], & \text{ if }|I_k|=2s\text{ for some }s\in\mathbb{N}, \\
    [\underbrace{2,2,\dotsc,2}_{s+1\text{ terms}},\underbrace{\tfrac{1}{2},\tfrac{1}{2},\dotsc,\tfrac{1}{2}}_{s\text{ terms}}], & \text{ if }|I_k|=2s+1\text{ for some }s\in\mathbb{N}_0, 
    \end{cases}
\end{align*}
and 
\begin{align*}
\mathcal{J}_{k}=
\begin{cases}
[\tfrac{1}{4},1,\dotsc,1], & \text{if }  k=1,\\
[\tfrac{1}{2},1,\dotsc,1], &\text{if else},
\end{cases}
\end{align*}
with $|\mathcal{J}_{k}|=|J_k|$. 
We consider the weighted backward shift $B_w\colon X\to X$ defined by the weight 
\[
w= (\underbrace{\dotsc\tfrac{1}{2},\tfrac{1}{2},\tfrac{1}{2}}_{-\mathbb{N}},\underbrace{\mathcal{I}_1,\mathcal{J}_1,\dotsc,\mathcal{I}_t, \mathcal{J}_t,2,2,\dotsc}_{\mathbb{N}_0}). 
\]
Since both $(w_i)_{i\in\mathbb{Z}}$ and $(\tfrac{1}{w_{i+1}})_{i\in\mathbb{Z}}$ are bounded, $B_w$ is an invertible operator on $X$. 
Clearly $B_w$ is strongly mixing (see \cite[Example 4.15]{G11}) and thus $\mathcal F$-mixing. 

\textbf{Case B: }$F\in \widetilde{\mathcal{F}}\setminus \mathcal{F}_{\mathrm{cf}}$. 
Write $F$ as the disjoint union of a sequence $(I_k)_{k\in\mathbb{N}}$ of maximal finite integer intervals,
and write $\mathbb{N}_0\setminus F$ as the disjoint union of a sequence $(J_k)_{k\in\mathbb{N}}$ of maximal finite integer intervals.
Without loss of generality we assume that the sequences $(\min I_k)_{k\in\mathbb{N}}$ and $(\min J_k)_{k\in\mathbb{N}}$ are increasing. 
Note that $0\in I_1$. 
For each $k\in\mathbb{N}$, let $\mathcal{I}_k$ and $\mathcal{J}_{k}$ be defined by the same formulas as in Case A.
Let $B_w$ be the bilateral weighted backward shift on $X$ defined by the weight
\[
    w=(\underbrace{\dotsc\tfrac{1}{2},\tfrac{1}{2},\tfrac{1}{2}}_{-\mathbb{N}},\underbrace{\mathcal{I}_1,\mathcal{J}_1,\dotsc,\mathcal{I}_n, \mathcal{J}_n,\mathcal{I}_{n+1}, \mathcal{J}_{n+1},\dotsc}_{\mathbb{N}_0}). 
\]
Since both $(w_i)_{i\in\mathbb{Z}}$ and $(\tfrac{1}{w_{i+1}})_{i\in\mathbb{Z}}$ are bounded, $B_w$ is an invertible operator on $X$. 
We first prove the following claim.

\medskip 

\noindent\textbf{Claim: }
If $m\in\mathbb{N}_0$ satisfies $\{m-r-1,m-r,\dotsc,m+r+1\}\subset F$ for some $r\in\mathbb{N}_0$ with $r>|I_1|$, then $w_0w_1\dotsb w_m\geq 2^r$. 

\begin{proof}[Proof of the Claim]
Denote $\varphi(j)=w_0w_1\dotsb w_j$ for $j\in\mathbb{N}_0$.
Since $(I_k)_{k\in\mathbb{N}}$ are disjoint, there is a unique $k\in\mathbb{N}$ such that $\{m-r-1,m-r,\dotsc,m+r+1\}\subset I_k$. 
Denote $a=\min I_k$ and $b=\max I_k$.
Since $|I_k|\geq 2r+3>|I_1|$, we have $k\geq 2$ and $|I_k|\geq 3$. 
If $|I_k|=2s$ for some $s\in\mathbb{N}$,
then 
\[
\varphi(a+i)=\begin{cases}
2^{i}, &0\leq i\leq s-1,\\
2^{2s-i-2}, &s\leq i\leq 2s-2;\\
1, &i=2s-1.
\end{cases}
\]
If $|I_k|=2s+1$ for some $s\in\mathbb{N}$,
then 
\[
\varphi(a+i)=\begin{cases}
2^{i}, &0\leq i\leq s,\\
2^{2s-i}, &s+1\leq i\leq 2s.
\end{cases}
\]
As $a+r+1\leq m \leq b-r-1$, it is clear that $\varphi(m)\geq 2^r$.
\end{proof}

To show that $B_w$ is $\mathcal{F}$-mixing, by Proposition \ref{Bes-char} it suffices to show that for any $M>0$ and $j\in\mathbb{N}$, $A_{M,j}\cap \overline{A}_{M,j}\in\mathcal{F}$, where 
\[
    A_{M,j}=\{n\in\mathbb{N}\colon  w_{j+1}w_{j+2}\dotsb w_{j+n}>M\}, \overline{A}_{M,j}=\{n\in\mathbb{N}\colon w_{j}w_{j-1}\dotsb w_{j-n+1}<\tfrac{1}{M}\}.
\]
Now fix $M>0$ and $j\in\mathbb{N}$. 
There is $q\in\mathbb{N}_0$ with $q>|I_1|$ and $2^q>Mw_0w_1w_2\dotsb w_j$. 
Consider the set 
\[
    E=\bigl\{m\in\mathbb{N}_0\colon \{m-q-j-1,m-q-j,\dotsc, m+q+j+1\}\subset F\bigr\}. 
\]
Recall that $\widetilde{\mathcal{F}}$ is the maximal shift-invariant thick subfamily of $\mathcal{F}$ and $F\in \widetilde{\mathcal{F}}$, by Remark \ref{thick}, $E\in\widetilde{\mathcal F}\subset \mathcal{F}$. 
Fix $n\in E$. 
We have $\{(j+n)-q-1,(j+n)-q,\dotsc, (j+n)+q+1\}\subset F$. 
By the Claim, 
\[
w_0w_1w_2\cdots w_{j+n}\geq 2^q. 
\]
This implies that 
\[
w_{j+1}w_{j+2}\cdots w_{j+n}=\frac{w_0w_1w_2\cdots w_{j+n}}{w_0w_1w_2\cdots w_{j}}\geq \frac{2^q}{w_0w_1w_2\cdots w_{j}}>M. 
\]
It is obvious that $n\geq q+j+1$. 
Then 
\begin{align*}
w_jw_{j-1}\dotsb w_{j-n+1}&=(w_jw_{j-1}\dotsb w_0)\cdot (w_{-1}\dotsb w_{j-n+1})\\
&=(w_0w_1\dotsb w_j)\cdot \frac{1}{2^{n-j-1}}  \leq  (w_0w_1\dotsb w_j)\cdot \frac{1}{2^{q}}<\frac{1}{M}.
\end{align*}
So $n\in A_{M,j}\cap \overline{A}_{M,j}$ and then $E\subset A_{M,j}\cap \overline{A}_{M,j}$. 
This gives that $A_{M,j}\cap \overline{A}_{M,j}\in\mathcal{F}$. 
Hence $B_w$ is $\mathcal{F}$-mixing, by Proposition \ref{Bes-char}. 

\medskip 

So in both cases, we have constructed an $\mathcal{F}$-mixing bilateral weighted backward shift $B_w$ on $X$, and for each $n\in\mathbb N_0$ one has
\begin{enumerate}
    \item $w_0w_1w_2\dotsb w_n =\tfrac{1}{2}$ if and only if $n\notin F$,
    \item $w_0w_1w_2\dotsb w_n \geq 1$ if and only if $n\in F$.  
\end{enumerate}
For each $x=\sum_{i\in\mathbb Z}x_ie_i\in X$, let 
\[
  \|x\|_+=\biggl\|\sum_{i\in\mathbb N_0}x_ie_i\biggr\|.   
\]
The function $\|\cdot\|_+$ is clearly a continuous seminorm on $X$. 
Consider the following neighborhood $W$ of $e_0$: 
\[
W=\bigl\{y\in X\colon \|y-e_0\|_+<\tfrac{4}{5}\bigr\}. 
\]

Now we verify that $N_{B_w}(W,W)=F$. 
On the one hand, fix $m\in N_{B_w}(W,W)\setminus\{0\}$. 
Then there is a vector $x=\sum_{i\in\mathbb Z}x_ie_i\in W$
such that $B_w^mx\in W$, which gives that $|x_m|<\frac{4}{5}$ and $|w_1w_2\cdots w_mx_m-1|<\frac{4}{5}$. 
Then we have $w_1w_2\cdots w_m|x_m|>\tfrac{1}{5}$ and 
\[
w_0w_1w_2\cdots w_m=\frac{w_0w_1w_2\cdots w_m|x_m|}{|x_m|}>\frac{w_0}{4}=\frac{1}{2}, 
\]
which gives that $m\in F$. 
On the other hand, fix $m\in F\setminus\{0\}$. 
Let $x=e_0+\tfrac{3}{10w_1w_2\dotsb w_m}e_m$. 
Then 
\[
\|x-e_0\|_+=\frac{3}{10w_1w_2\dotsb w_m}\|e_m\|=\frac{3w_0}{10w_0w_1w_2\dotsb w_m}\leq \frac{3w_0}{10}=\frac{3}{5}<\frac{4}{5}, 
\]
and hence $x\in W$. 
Moreover,
\[
\|B_w^mx-e_0\|_+=\Bigl\|\frac{3}{10}e_0-e_0 + (w_{-m+1}w_{-m+2}\dotsb w_{0}) e_{-m}\Bigr\|_+=\frac{7}{10}\|e_0\|=\frac{7}{10}<\frac{4}{5}. 
\]
So $B_w^mx\in W$. Therefore, $m\in N_{B_w}(W,W)$. 
\end{proof}

\begin{rem}
Given a full family $\mathcal{F}$ and $F\in \widetilde{\mathcal{F}}$, by a strategy similar to that used in the proof of Lemma \ref{L-WAG} one can find an $\mathcal{F}$-mixing unilateral weighted backward shift $B_w$ on $X$, 
and a non-empty open subset $W$ of $X$ such that $N_{B_w}(W,W)= F\cup\{0\}$, where $X=c_0(\mathbb{N}_0)$ or $ \ell^p(\mathbb{N}_0)$, $1\leq p<\infty$. 
\end{rem}

By Lemma \ref{L-WAG}, for a given full family $\mathcal{F}$, we characterize the operators which are weakly disjoint from all $\mathcal{F}$-mixing operators, which is an analogue of the Weiss-Akin-Glasner Theorem (Theorem~\ref{thm:WAG-revised}) from topological dynamics.
    
\begin{thm}\label{linearWAG}
Let $T\colon X\to X$ be an operator and let $\mathcal{F}$ be a full family. 
Then the following assertions are equivalent: 
    \begin{enumerate}
        \item $T$ is $\bigl(\widetilde{\mathcal{F}}\bigr)^*$-transitive; 
        \item $T$ is weakly disjoint from all $\mathcal{F}$-mixing operators; 
        \item $T$ is weakly disjoint from all invertible $\mathcal{F}$-mixing bilateral weighted backward shifts on $Y$,
        where
        $Y=c_0(\mathbb{Z})$ or $\ell^p(\mathbb{Z}),1\leq p<\infty$.
   \end{enumerate}
\end{thm}

\begin{proof}
The implications (1)$\Rightarrow$(2)$\Rightarrow$(3) are clear. 

(3)$\Rightarrow$(1). Let $U, V$ be a pair of non-empty open subsets of $X$. 
Our goal is to show that $N_T(U,V)$ intersects each member in $\widetilde{\mathcal{F}}$. 
Fix $F\in\widetilde{\mathcal{F}}$. By Lemma~\ref{L-WAG}, there exists an invertible $\mathcal{F}$-mixing bilateral weighted backward shift $S$ on $Y$,  and a non-empty open subset $W$ of $Y$ such that $N_S(W,W)= F\cup\{0\}$. 
 Since $T\oplus S$ is hypercyclic, $N_T(U,V)\cap N_S(W,W)$ is an infinite set, while $N_T(U,V)\cap N_S(W,W)= N_T(U,V)\cap (F\cup\{0\})=(N_T(U,V)\cap F)\cup (N_T(U,V)\cap\{0\})$. Consequently, $N_T(U,V)\cap F$ is an infinite set, and then $N_T(U,V)\in\bigl(\widetilde{\mathcal{F}}\bigr)^*$. 
Hence $T$ is $\bigl(\widetilde{\mathcal{F}}\bigr)^*$-transitive. 
\end{proof}

For a shift-invariant thick family $\mathcal{F}$,
if every $\mathcal{F}^*$-transitive operator is weakly mixing, then the characterization of the weak disjointness of $\mathcal{F}$-transitive operators can be strengthened as follows. 

\begin{prop}\label{prop:F-star-weak-mixing}
Let $T\colon X\to X$ be an operator and let $\mathcal{F}$ be a shift-invariant thick family. 
If every $\mathcal{F}^*$-transitive operator is weakly mixing, then the following assertions are equivalent: 
\begin{enumerate}
    \item $T$ is $\mathcal{F}^*$-transitive;
    \item for any $\mathcal{F}$-transitive operator $S\colon Y\to Y$, $T\oplus S$ is weakly mixing; 
    \item for any invertible $\mathcal{F}$-transitive bilateral weighted backward shift $S$ on $Y$, 
    where $Y=c_0(\mathbb{Z})$, or $\ell^p(\mathbb{Z})$, $1\leq p<\infty$, $T\oplus S$ is weakly mixing. 
\end{enumerate}
\end{prop}
\begin{proof}
    The implication (2)$\Rightarrow$(3) is clear. (3)$\Rightarrow$(1) follows from Theorem \ref{linearWAG}. 
    
     (1)$\Rightarrow$(2). Let $S\colon Y\to Y$ be an $\mathcal{F}$-transitive operator. Fix non-empty open subsets $U_1, V_1$ of $X$, and $U_2, V_2$ of $Y$. Let $k\in\mathbb{N}_0$. 
     By the hypothesis, $T$ is weakly mixing and $\mathcal{F}^*$-transitive.
    So $T$ is $\mathcal{F}^*$-mixing. 
     Then $A=\bigl\{n\in\mathbb{N}_0: \{n,n+1,\dotsc,n+k\}\subset N_T(U_1, V_1)\bigr\}\in\mathcal{F}^*$. 
     By the $\mathcal{F}$-transitivity of $S$ and the thickness of $\mathcal{F}$, we have $B=\bigl\{n\in\mathbb{N}_0: \{n,n+1,\dotsc,n+k\}\subset N_S(U_2, V_2)\bigr\}\in\mathcal{F}$. 
     Now $\emptyset\neq A\cap B\subset \{n\in\mathbb{N}_0: \{n,n+1,\dotsc,n+k\}\subset N_T(U_1,V_1)\cap N_S(U_2, V_2)\}$ and then $N_T(U_1,V_1)\cap N_S(U_2, V_2)\in\mathcal{F}_{\mathrm{t}}$. 
     By Proposition \ref{F-intersection}, $T\oplus S$ is weakly mixing. 
\end{proof}

Since $\mathcal{F}_{\mathrm{ts}}$ and $\mathcal{F}_{\mathrm{t}}$ are shift-invariant thick families, and topologically ergodic operators are weakly mixing, by Theorem \ref{linearWAG} and Proposition~\ref{prop:F-star-weak-mixing}, we can strengthen Cardeccia's results (Theorems \ref{wk-mixing-dsj} and \ref{ps-dsj}) as follows. 

\begin{coro}\label{cor:top-ergodic}
Let $T\colon X\to X$ be an operator.
Then the following assertions are equivalent: 
   \begin{enumerate}
       \item $T$ is topologically ergodic; 
       \item $T$ is weakly disjoint from all weakly mixing operators; 
       \item $T$ is weakly disjoint from  all invertible weakly mixing bilateral weighted backward shifts on $Y$,
        where $Y=c_0(\mathbb{Z})$ or $\ell^p(\mathbb{Z}),1\leq p<\infty$; 
       \item for any invertible weakly mixing bilateral weighted backward shift $S$ on $Y$, $T\oplus S$ is weakly mixing, where $Y=c_0(\mathbb{Z})$ or $\ell^p(\mathbb{Z}),1\leq p<\infty$; 
       \item for any weakly mixing operator $S\colon Y\to Y$, $T\oplus S$ is weakly mixing, where $Y=c_0(\mathbb{Z})$ or $\ell^p(\mathbb{Z}),1\leq p<\infty$. 
   \end{enumerate} 
\end{coro}

\begin{coro}\label{cor:ps-trans}
Let $T\colon X\to X$ be an operator.
 Then the following assertions are equivalent: 
   \begin{enumerate}
       \item $T$ is $\mathcal{F}_{\mathrm{ps}}$-transitive; 
       \item $T$ is weakly disjoint from all topologically ergodic operators; 
       \item $T$ is weakly disjoint from all invertible topologically ergodic bilateral weighted backward shifts on $Y$,
        where $Y=c_0(\mathbb{Z})$ or $\ell^p(\mathbb{Z}),1\leq p<\infty$.
   \end{enumerate} 
\end{coro}

As we can see, if each member of a family $\mathcal{F}$ contains a return time set of an $\mathcal{F}$-transitive operator, then the $\mathcal{F}^*$-transitivity of an operator is equivalent to the weak disjointness from all $\mathcal{F}$-transitive operators. 
By Lemma \ref{L-WAG}, every member of a shift-invariant thick family contains a return time set of an $\mathcal{F}$-mixing operator. 
However, one can construct a family $\mathcal{F}$ where every member contains a return time set of an $\mathcal{F}$-transitive operator but is not thick, using the following result. 
    
\begin{thm}[{\cite[Theorem 4.1]{DR09} or \cite[Theorem 1.2]{BE09}}]
There is a hypercyclic operator which is not weakly mixing. 
\end{thm}
  
Let $T\colon X\to X$ be a hypercyclic operator which is not weakly mixing and let $\mathcal{N}_T=\{A\subset\mathbb{N}_0\colon N(U,V)\subset A\text{ for some non-empty open subsets }U,V\text{ of }X\}$. 
It is easy to check that an operator $S\colon Y\to Y$ is $\mathcal{N}_T^*$-transitive if and only if it is weakly disjoint from all $\mathcal{N}_T$-transitive operators if and only if it is weakly disjoint from $T$.
But since $T$ is not weakly mixing, we have $N(U_0,V_0)\not\in\mathcal{F}_{\mathrm{t}}$ for some non-empty open subsets $U_0,V_0$ of $X$ and then $\mathcal{N}_T$ is not a thick family. 
Lemma~\ref{L-WAG} gives a characterization of subsets of non-negative integers which can be realized by the return time sets of $\mathcal{F}$-mixing operators.
But the characterization of return time sets remains open for hypercyclic operators which are not weakly mixing.
This raises the following question: 

\begin{Ques}
    Besides shift-invariant thick families, what kind of family $\mathcal{F}$ have the property that every member of $\mathcal{F}$ contains a return time set of an $\mathcal{F}$-transitive operator? 
\end{Ques}

By the Weiss-Akin-Glasner Theorem (Theorem \ref{WAG-thm}), there is no family $\mathcal{F}$ such that the $\mathcal{F}$-transitivity of compact dynamical systems is equivalent to the weak disjointness from all $\mathcal{F}$-transitive compact dynamical systems, since there exists a topologically ergodic compact dynamical system which is not weakly mixing. 
But things change under the framework of linear dynamical systems, since each topologically ergodic operator must be $\mathcal{F}_{\mathrm{ts}}$-transitive, thus weakly mixing, so we ask: 

\begin{Ques}
    Is there a family $\mathcal{F}$ such that 
    an operator is $\mathcal{F}$-transitive 
    if and only if it is weakly disjoint from all $\mathcal{F}$-transitive operators? 
\end{Ques}

\begin{rem}
It is clear that every weakly mixing operator is $\mathcal{F}_{\mathrm{ps}}$-transitive.
A natural question is whether these two classes of operators are distinct, see \cite[Question in page 1328]{C24}. 
We note that by Corollaries \ref{cor:top-ergodic} and \ref{cor:ps-trans}, 
an operator is weakly disjoint from all weakly mixing operators  if and only if it is weakly disjoint from all $\mathcal{F}_{\mathrm{ps}}$-transitive operators.
\end{rem}

In 1995, Salas provided in \cite[Corollary 2.3]{S95} a hypercyclic operator $T$ on $\ell^2(\mathbb{Z})$ such that $T^*$ is also hypercyclic. 
In Lemma \ref{L-WAG}, if we restrict the underlying space $X$ to be $\ell^2(\mathbb{Z})$, 
require that $\mathbb{N}_0\setminus F$ is large enough, 
and modify our construction slightly, 
then the operator we obtain possesses a hypercyclic adjoint. 

\begin{prop}\label{adjoint}
Let $\mathcal{F}_1$ and $\mathcal{F}_2$ be shift-invariant thick families. If there are $F_1\in\mathcal{F}_1,F_2\in\mathcal{F}_2$ with $F_1\cap F_2=\emptyset$, then there is an invertible operator $T$ on $\ell^2(\mathbb{Z})$ such that 
$T$ is $\mathcal{F}_1$-transitive and 
$T^*$ is $\mathcal{F}_2$-transitive.  
\end{prop}

\begin{proof}
Let $F=F_1$. It is clear that both $F$ and $\mathbb{N}_0\setminus F$ are not cofinite. 
Write $F$ as the disjoint union of a sequence $(I_k)_{k\in\mathbb{N}}$ of maximal finite integer intervals. 
Since $\bigl\{ n\in\mathbb{N}_0\colon \{n-2,n-1,n\}\subset F\bigr\}\in\mathcal{F}_1$, without loss of generality, we require that $0\notin F,1\notin F$ and $|I_k|\geq 3$ for each $k\in\mathbb{N}$. 
Also we assume that the sequence $(\operatorname{min}I_k)_{k\in\mathbb{N}}$ is increasing. 
Let $J_1=[1,\operatorname{min}I_1-1]\cap\mathbb{N}_0$. 
For each $k\in\mathbb{N}$, let $J_{k+1}=[\operatorname{max}I_{k}+1,\operatorname{min}I_{k+1}-1]\cap\mathbb{N}_0$. 
For $n\in\mathbb{N}$, let 
\begin{align*}
    \mathcal{I}_n=\begin{cases}
        [\underbrace{2,2,\dotsc,2}_{m\text{ terms}},\underbrace{\tfrac{1}{2},\tfrac{1}{2},\dotsc,\tfrac{1}{2}}_{m\text{ terms}}], \text{ if }|I_n|=2m\text{ for some }m\in\mathbb{N}, \\
        [\underbrace{2,2,\dotsc,2}_{m\text{ terms}},\underbrace{\tfrac{1}{2},\tfrac{1}{2},\dotsc,\tfrac{1}{2}}_{m\text{ terms}},1], \text{ if }|I_n|=2m+1\text{ for some }m\in\mathbb{N}, 
    \end{cases}
\end{align*}
and 
\begin{align*}
    \mathcal{J}_n=\begin{cases}
        [\underbrace{\tfrac{1}{2},\tfrac{1}{2},\dotsc,\tfrac{1}{2}}_{m\text{ terms}},\underbrace{2,2,\dotsc,2}_{m\text{ terms}}], \text{ if }|J_n|=2m\text{ for some }m\in\mathbb{N}, \\
        [\underbrace{\tfrac{1}{2},\tfrac{1}{2},\dotsc,\tfrac{1}{2}}_{m\text{ terms}},\underbrace{2,2,\dotsc,2}_{m\text{ terms}},1], \text{ if }|J_n|=2m+1\text{ for some }m\in\mathbb{N}, \\
        [1], \text{ if } |J_n|=1. 
    \end{cases}
\end{align*}
Also, let 
\begin{align*}
    \overline{\mathcal{I}_n}=\begin{cases}
        [\underbrace{2,2,\dotsc,2}_{m\text{ terms}},\underbrace{\tfrac{1}{2},\tfrac{1}{2},\dotsc,\tfrac{1}{2}}_{m\text{ terms}}], \text{ if }|I_n|=2m\text{ for some }m\in\mathbb{N}, \\
        [1,\underbrace{2,2,\dotsc,2}_{m\text{ terms}},\underbrace{\tfrac{1}{2},\tfrac{1}{2},\dotsc,\tfrac{1}{2}}_{m\text{ terms}}], \text{ if }|I_n|=2m+1\text{ for some }m\in\mathbb{N}, 
    \end{cases}
\end{align*}
and 
\begin{align*}
    \overline{\mathcal{J}_n}=\begin{cases}
        [\underbrace{\tfrac{1}{2},\tfrac{1}{2},\dotsc,\tfrac{1}{2}}_{m\text{ terms}},\underbrace{2,2,\dotsc,2}_{m\text{ terms}}], \text{ if }|J_n|=2m\text{ for some }m\in\mathbb{N}, \\
        [1,\underbrace{\tfrac{1}{2},\tfrac{1}{2},\dotsc,\tfrac{1}{2}}_{m\text{ terms}},\underbrace{2,2,\dotsc,2}_{m\text{ terms}}], \text{ if }|J_n|=2m+1\text{ for some }m\in\mathbb{N}, \\
        [1], \text{ if } |J_n|=1. 
    \end{cases}
\end{align*}
Let $B_w$ be the bilateral weighted backward shift on $X$ defined by the weight $w$, where $w=(w_i)_{i\in\mathbb{Z}}$ is given by 
\[
    w=(\underbrace{\dotsc \overline{\mathcal{I}_n},\overline{\mathcal{J}_n},\dotsc\overline{\mathcal{I}_1},\overline{\mathcal{J}_1}}_{-\mathbb{N}},
    \underbrace{[1]}_{0},
    \underbrace{\mathcal{J}_1,\mathcal{I}_1,\dotsc, \mathcal{J}_n,\mathcal{I}_n,\dotsc}_{\mathbb{N}}). 
\]
For each $i\in -\mathbb{N}$, we have $w_i=\tfrac{1}{w_{-i}}$. 
Since both $(w_i)_{i\in\mathbb{Z}}$ and $(\tfrac{1}{w_{i+1}})_{i\in\mathbb{Z}}$ are bounded, $B_w$ is an invertible operator on $\ell^2(\mathbb{Z})$. 
Let $v=(w_{1-i})_{i\in\mathbb{N}_0}$ and $T=B_w$. 
It is easy to see that $T^*$ is conjugate to $B_v$. 
It remains to show that $B_w$ is $\mathcal{F}_1$-transitive and $B_v$ is $\mathcal{F}_2$-transitive. 

By the same argument as the one in Case B of the proof of Lemma \ref{L-WAG}, we have the following facts: 
\begin{itemize}
    \item if $m\in\mathbb{N}_0$ satisfies $\{m-r,m-r+1,\dotsc,m+r+1\}\subset F$ for some $r\in\mathbb{N}_0$, then $w_1\dotsb w_m\geq 2^r$; 
    \item if $m\in\mathbb{N}_0$ satisfies $\{m-r,m-r+1,\dotsc,m+r+1\}\subset \mathbb{N}\setminus F$ for some $r\in\mathbb{N}_0$, then $w_1\dotsb w_m\leq \tfrac{1}{2^r}$. 
\end{itemize}

To show that $B_w$ is $\mathcal{F}_1$-transitive,
by Proposition~\ref{Bes-char} it suffices to show that for any $M>0$ and $j\in\mathbb{N}$, $A_{M,j}(w)\cap \overline{A}_{M,j}(w)\in\mathcal{F}_1$, where 
        \[
            A_{M,j}(w)=\{n\in\mathbb{N}: w_{j+1}w_{j+2}\dotsb w_{j+n}>M\}, \overline{A}_{M,j}(w)=\{n\in\mathbb{N}: w_{j}w_{j-1}\dotsb w_{j-n+1}<\tfrac{1}{M}\}.
        \]
Fix $j\in\mathbb{N}$ and $M>0$. 
There is $q_1\in\mathbb{N}$ with $2^{q_1}>Mw_1\dotsb w_j$. 
Consider the set  
\[
    E_1=\bigl\{m\in\mathbb{N}_0: \{m-q_1-j-1, m-q_1-j,\dotsc,m+q_1+j+1\}\subset F\bigr\}. 
\]
Since $\mathcal{F}_1$ is a shift-invariant thick family and $F\in\mathcal{F}_1$, by Remark \ref{thick}, $E_1\in\mathcal{F}_1$. 
Fix $n\in E_1$, we have $\{(n+j)-q_1, (n+j)-q_1+1,\dotsc,(n+j)+q_1+1\}\subset F$ and $\{(n-j-1)-q_1, (n-j-1)-q_1+1,\dotsc,(n-j-1)+q_1+1\}\subset F$. 
Then $w_1w_2\dotsb w_{j+n}\geq 2^{q_1}$ and $w_1w_2\dotsb w_{n-j-1}\geq 2^{q_1}$. 
We have 
\[
    w_{j+1}w_{j+2}\dotsb w_{j+n}=\frac{w_1w_2\dotsb w_{j+n}}{w_1w_2\dotsb w_j}\geq \frac{2^{q_1}}{w_1w_2\dotsb w_j}>M, 
\]
and 
\[
    w_jw_{j-1}\dotsb w_0w_{-1}\dotsb w_{j-n+1}=w_0w_1\dotsb w_j \frac{1}{w_{1}\dotsb w_{n-j-1}}\leq \frac{w_1\dotsb w_j}{2^{q_1}}<\frac{1}{M}. 
\]
This gives that $n\in A_{M,j}(w)\cap \overline{A}_{M,j}(w)$ and then $E_1\subset A_{M,j}(w)\cap \overline{A}_{M,j}(w)$.  
Then $A_{M,j}(w)\cap \overline{A}_{M,j}(w)\in\mathcal{F}_1$, and $B_w$ is $\mathcal{F}_1$-transitive by Proposition \ref{Bes-char}. 

To show that $B_v$ is $\mathcal{F}_2$-transitive,
by Proposition~\ref{Bes-char} it suffices to show that for any $M>0$ and $j\in\mathbb{N}$, $A_{M,j}(v)\cap \overline{A}_{M,j}(v)\in\mathcal{F}_2$, where 
        \[
            A_{M,j}(v)=\{n\in\mathbb{N}: v_{j+1}v_{j+2}\dotsb v_{j+n}>M\}, \overline{A}_{M,j}(v)=\{n\in\mathbb{N}: v_{j}v_{j-1}\dotsb v_{j-n+1}<\tfrac{1}{M}\}.
        \]
Fix $j\in\mathbb{N}$ and $M>0$. 
There is $q_2\in\mathbb{N}$ with $2^{q_2}>M\operatorname{max}\{w_0w_{-1}\dotsb w_{-j+1},v_jv_{j-1}\dotsb v_2\}$. 
Consider the set  
\[
    E_2=\bigl\{m\in\mathbb{N}_0: \{m-q_2-j,m-q_2-j+1,\dotsc, m+q_2+j+1\}\subset \mathbb{N}\setminus F\bigr\}. 
\]
Since $\mathbb{N}_0\setminus F\in\mathcal{F}_2$ and $\mathcal{F}_2$ is full, we have $\mathbb{N}\setminus F\in\mathcal{F}_2$. 
Since $\mathcal{F}_2$ is a shift-invariant thick family, by Remark \ref{thick}, $E_2\in\mathcal{F}_2$. 
Fix $n\in E_2$.
We have $\{(n+j)-q_2,(n+j)-q_2+1,\dotsc, (n+j)+q_2+1\}\subset \mathbb{N}\setminus F$ and $\{(n-j)-q_2,(n-j)-q_2+1,\dotsc, (n-j)+q_2+1\}\subset \mathbb{N}\setminus F$. 
By the above analysis, we have $w_1w_2\dotsb w_{j+n-1}\leq\frac{1}{2^{q_2}}$ and $w_1w_2\dotsb w_{n-j}\leq\frac{1}{2^{q_2}}$. Note that 
\[
    v_1\dotsb v_{j+n}=\frac{1}{w_1w_2\dotsb w_{j+n-1}}=\frac{1}{w_0w_1w_2\dotsb w_{j+n-1}}. 
\]
We have 
\[
    v_{j+1}\dotsb v_{j+n}=\frac{\frac{1}{w_1w_2\dotsb w_{j+n-1}}}{v_1\dotsb v_j}\geq\frac{2^{q_2}}{w_0w_{-1}\dotsb w_{-j+1}}>M, 
\]
and 
\[
    v_jv_{j-1}\dotsb v_0 v_{-1}\dotsb v_{j-n+1}=v_jv_{j-1}\dotsb v_2 w_0w_1 w_2\dotsb w_{-j+n}\leq\frac{v_jv_{j-1}\dotsb v_2}{2^{q_2}}<\frac{1}{M}. 
\]
This gives that $n\in A_{M,j}(v)\cap \overline{A}_{M,j}(v)$ and then $E_2\subset A_{M,j}(v)\cap \overline{A}_{M,j}(v)$. 
Then $A_{M,j}(v)\cap \overline{A}_{M,j}(v)\in\mathcal{F}_2$, and $B_v$ is $\mathcal{F}_2$-transitive by Proposition \ref{Bes-char}. 
\end{proof}

\subsection{Some Examples}

Applying Lemma \ref{L-WAG}, we can provide a large number of operators with some prescribed properties. 
We first give two general results and then provide certain examples. 

\begin{prop}\label{distinguish}
Let $\mathcal{F}_1$ be a shift-invariant thick family and $\mathcal{F}_2$ be a family.
Then there exists an invertible $\mathcal{F}_1$-transitive operator on $X$ which is not $\mathcal{F}_2$-transitive if and only if $\mathcal{F}_1\setminus \mathcal{F}_2\neq\emptyset$, where $X=c_0(\mathbb{Z})$ or $\ell^p(\mathbb{Z})$, $1\leq p<\infty$. 
\end{prop}

\begin{proof}
On the one hand, let $T\colon X\to X$ be an invertible $\mathcal{F}_1$-transitive operator which is not $\mathcal{F}_2$-transitive. 
Then there are non-empty open subsets $U,V$ of $X$ with $N(U,V)\notin \mathcal{F}_2$. 
We are done since $N(U,V)\in \mathcal{F}_1$. 

Conversely, fix $F\in \mathcal{F}_1\setminus \mathcal{F}_2$. 
By Lemma \ref{L-WAG}, there exist an invertible $\mathcal{F}_1$-transitive operator $T\colon X\to X$ and a non-empty open subset $W$ of $X$ such that 
\[
    N(W,W)= F\cup\{0\}, 
\] 
where $X=c_0(\mathbb{Z})$ or $\ell^p(\mathbb{Z})$, $1\leq p<\infty$. 
Since $X$ has no isolated points, there are non-empty open sets $W', W''\subset W$ with $N(W', W'')\subset [1,+\infty)$. 
Finally we have $N(W', W'')\subset F$. 
Since $F\notin \mathcal F_2$, we have that $N(W', W'')\notin\mathcal F_2$. 
Hence $T$ is not $\mathcal{F}_2$-transitive. 
\end{proof}

\begin{prop}\label{construction}
    Let $\mathcal{F}_1$ and $\mathcal{F}_2$ be shift-invariant thick families.
    Then there exists an invertible $\mathcal{F}_1$-transitive operator $T\colon X\to X$ and an invertible $\mathcal{F}_2$-transitive operator $S\colon Y\to Y$ 
    such that $T$ is not weakly disjoint from $S$ if and only if there are $F_1\in\mathcal{F}_1$ and $F_2\in\mathcal{F}_2$ such that $F_1\cap F_2=\emptyset$.
\end{prop}

\begin{proof}    
On the one hand, let $T\colon X\to X$ be an invertible $\mathcal{F}_1$-transitive operator and $S\colon Y\to Y$ be an invertible $\mathcal{F}_2$-transitive operator such that $T$ is not weakly disjoint from $S$. Then there are non-empty open subsets $U_1,V_1$ of $X$  and $U_2,V_2$ of $Y$ with $N_T(U_1,V_1)\cap N_S(U_2,V_2)=\emptyset$. We are done since $N_T(U_1,V_1)\in \mathcal{F}_1$ and $N_S(U_2,V_2)\in \mathcal{F}_2$. 

Conversely, by the hypothesis, $F_1\in\mathcal{F}_1\setminus \mathcal{F}^*_2$. By Proposition \ref{distinguish}, there exists an invertible $\mathcal{F}_1$-transitive operator which is not $\mathcal{F}^*_2$-transitive. Since $\mathcal{F}_2$ is a shift-invariant thick family, by Theorem \ref{linearWAG}, there exists an invertible $\mathcal{F}_2$-transitive operator $S\colon Y\to Y$ such that $T$ is not weakly disjoint from $S$. 
\end{proof}

In 2019, B\`es et al.\@ have constructed several operators to distinguish transitivity with respect to  different families \cite{BMPP19}. 
We provide here more invertible operators of this kind.  
We recall the definitions of density and Banach density of a subset of $\mathbb{N}_0$. 

\begin{defn}
Let $A\subset\mathbb{N}_0$. 
\begin{itemize}
    \item The \emph{upper and lower density} of $A$ are given by 
    \[
    \overline{\operatorname{D}}(A)=\limsup_{n\to \infty} \frac{|\{j\in A: j\leq n-1\}|}{n},\ \underline{\operatorname{D}}(A)=\liminf_{n\to \infty} \frac{|\{j\in A: j\leq n-1\}|}{n};   
    \]
    \item The \emph{upper and lower Banach density} of $A$ are given by \[
    \overline{\operatorname{BD}}(A)=\limsup_{n-m\to \infty} \frac{|\{j\in A: m\leq j\leq n-1\}|}{n-m}, 
    \ \underline{\operatorname{BD}}(A)=\liminf_{n-m\to \infty} \frac{|\{j\in A: m\leq j\leq n-1\}|}{n-m}. 
    \]
\end{itemize}
\end{defn}

It is observed that for any $A\subset\mathbb{N}_0$, we have $\underline{\operatorname{BD}}(A)\leq\underline{\operatorname{D}}(A)\leq\overline{\operatorname{D}}(A)\leq\overline{\operatorname{BD}}(A)$. 
We denote 
\begin{itemize}
    \item  $\mathcal{F}_{\mathrm{pud}}=\bigl\{A\subset\mathbb{N}_0: \overline{\operatorname{D}}(A)>0\bigr\}$;  \quad\quad $\mathcal{F}_{\mathrm{pld}}=\bigl\{A\subset\mathbb{N}_0: \underline{\operatorname{D}}(A)>0\bigr\}$; 
    \item  $\mathcal{F}_{\mathrm{ud1}}=\bigl\{A\subset\mathbb{N}_0: \overline{\operatorname{D}}(A)=1\bigr\}$;  \quad\quad $\mathcal{F}_{\mathrm{ld1}}=\bigl\{A\subset\mathbb{N}_0: \underline{\operatorname{D}}(A)=1\bigr\}$; 
    \item  $\mathcal{F}_{\mathrm{pubd}}=\bigl\{A\subset\mathbb{N}_0: \overline{\operatorname{BD}}(A)>0\bigr\}$;\quad $\mathcal{F}_{\mathrm{lbd1}}=\bigl\{A\subset\mathbb{N}_0: \underline{\operatorname{BD}}(A)=1\bigr\}$. 
\end{itemize}

It can be checked that the above six families are shift-invariant and $\mathcal{F}_{\mathrm{ld1}}, \mathcal{F}_{\mathrm{lbd1}}$ are filters. 
By Remark \ref{thick}(5), both $\mathcal{F}_{\mathrm{ld1}}$ and $\mathcal{F}_{\mathrm{lbd1}}$ are shift-invariant thick families. 
Besides, we have $\mathcal{F}^*_{\mathrm{ld1}}=\mathcal{F}_{\mathrm{pud}}$ and $\mathcal{F}^*_{\mathrm{lbd1}}=\mathcal{F}_{\mathrm{pubd}}$. 
By Proposition \ref{distinguish}, once we have found a subset of $\mathbb{N}_0$ with some prescribed properties, we will obtain an operator to distinguish certain classes of hypercyclic operators. 

\begin{prop}\label{distinguish-examples}
Let $X=c_0(\mathbb{Z})$ or $\ell^p(\mathbb{Z})$, $1\leq p<\infty$. 
\begin{enumerate}
    \item There exists an invertible weakly mixing operator on $X$ which is not $\mathcal{F}_{\mathrm{pud}}$-transitive. 
    \item There exists an invertible topologically ergodic operator on $X$ which is not $\mathcal{F}^*_{\mathrm{ip}}$-transitive. 
    \item There exists an invertible topologically ergodic operator on $X$ which is not $\mathcal{F}_{\mathrm{ud1}}$-transitive. 
    \item There exists an invertible $\mathcal{F}_{\mathrm{ld1}}$-transitive operator on $X$ which is not topologically ergodic. 
    \item There exists an invertible $\mathcal{F}_{\mathrm{lbd1}}$-transitive operator on $X$ which is not strongly mixing. 
    \item There exists an invertible $\mathcal{F}^*_{\mathrm{ip}}$-transitive operator on $X$ which is not $\Delta^*$-transitive. 
    \item There exists an invertible $\Delta^*_{\mathrm{ip}}$-transitive operator on $X$ which is not $\mathcal{F}^*_{\mathrm{ip}}$-transitive. 
\end{enumerate}
\end{prop}

\begin{proof}
(1) Let $F=\bigcup_{n\in\mathbb{N}} \{2^n,2^n+1,\dotsc,2^n+n\}$.
It is clear that $F\in\mathcal{F}_{\mathrm{t}}\setminus \mathcal{F}_{\mathrm{pud}}$. 
Since $\mathcal{F}_{\mathrm{t}}$ is a shift-invariant thick family, the result follows from Proposition \ref{distinguish}. 

(2) Let $F=\mathbb{N}_0\setminus \operatorname{FS}((10^{n})_{n\in\mathbb{N}})$. 
It can be verified that $F\in\mathcal{F}_{\mathrm{ts}}\setminus\mathcal{F}^*_{\mathrm{ip}}$. 
Since $\mathcal{F}_{\mathrm{ts}}$ is a shift-invariant thick family, the result follows from Proposition \ref{distinguish}. 

(3) For each $n\in\mathbb{N}_0$, let $F_n = \{(n 3^n)k+i\colon k\in\mathbb{N}, i\in\{0,1,\dotsc,n-1\}\}$.
Then for each $n\in\mathbb{N}$, $F_n$ is syndetic and $\overline{\operatorname{D}}(F_n)=\underline{\operatorname{D}}(F_n)=\frac{1}{3^n}$. 
Then $F: =\bigcup_{n=1}^\infty F_n$ is thickly syndetic and $\overline{\operatorname{D}}(F)\leq\sum_{n=1}^\infty \overline{\operatorname{D}}(F_n)=\frac{2}{3}$.
This shows that $F\in\mathcal{F}_{\mathrm{ts}}\setminus\mathcal{F}_{\mathrm{ud1}}$. 
Since $\mathcal{F}_{\mathrm{ts}}$ is a shift-invariant thick family, the result follows from Proposition \ref{distinguish}. 

(4) Let $F=\mathbb{N}_0\setminus \bigcup_{n\in\mathbb{N}} \{2^n,2^n+1,\dotsc,2^n+n\}$. 
By the proof of (1), $F\in\mathcal{F}_{\mathrm{ld1}}\setminus\mathcal{F}_{\mathrm{s}}$. 
Since $\mathcal{F}_{\mathrm{ld1}}$ is a shift-invariant thick family, the result follows from Proposition \ref{distinguish}. 

(5) Let $F=\mathbb{N}_0\setminus \{2^n\colon n\in\mathbb{N}\}$.
It is clear that $F\in\mathcal{F}_{\mathrm{lbd1}}\setminus\mathcal{F}_{\mathrm{cf}}$. 
Since $\mathcal{F}_{\mathrm{lbd1}}$ is a shift-invariant thick family, the result follows from Proposition \ref{distinguish}. 

(6) Let $(s_n)_{n\in\mathbb{N}}$ be a sequence of positive integers such that $s_{n+1}>4(\sum_{i=1}^n s_i+n)$ for each $n\in\mathbb{N}$. 
Put $S=\{s_n\colon n\in\mathbb{N}\}$ and $A=\mathbb{N}_0\setminus (S-S)$. 
The family 
\[
    \mathcal{F}_A=\biggl\{B\subset\mathbb{N}_0: \text{there exist }m_1,m_2,\dotsc,m_k\in\mathbb{Z}\text{ with }\bigcap_{i=1}^k (A+m_i)\biggr\}
\]
is a shift-invariant thick subfamily of $\mathcal{F}^*_{\mathrm{ip}}$ (see \cite[Example B]{HY04}). 
Now we have $A\in \mathcal{F}_A\setminus\Delta^*$. 
Then the result follows from Proposition \ref{distinguish}. 

(7) Let $(s_n)_{n\in\mathbb{N}}$ be the same sequence as the one in the proof of (6). 
Put $A=\mathbb{N}_0\setminus \operatorname{FS}((s_n)_{n\in\mathbb{N}})$. 
For each $k\in\mathbb{N}$, let $A_k=\bigcap_{i=-k}^k (A+i)$. 
The family 
\[
    \mathcal{F}_A=\biggl\{B\subset\mathbb{N}_0: \text{there exists }k\in\mathbb{N}\text{ with }A_k\subset B\biggr\}
\]
is a shift-invariant thick subfamily of $\Delta^*_{\mathrm{ip}}$ (see \cite[Example C]{HY04}). 
Now we have $A\in\mathcal{F}_A\setminus \mathcal{F}^*_{\mathrm{ip}}$. 
Then the result follows from Proposition \ref{distinguish}. 
\end{proof}

We provide examples of weakly mixing operators without hypercyclic direct sums. 
    
\begin{coro}\label{hyper-not-dsj}
There exist an invertible weakly mixing operator $T$ and an invertible $\mathcal{F}_{\mathrm{ld1}}$-transitive operator $S$ such that $T$ is not weakly disjoint from $S$.
\end{coro}
\begin{proof}
Since there are $F_1\in\mathcal{F}_{\mathrm{t}}$ and $F_2\in\mathcal{F}_{\mathrm{ld1}}$ with $F_1\cap F_2=\emptyset$ (see the proof of Proposition \ref{distinguish-examples}(4)), the result follows from Proposition \ref{construction}. 
\end{proof}

In Section 3 we have demonstrated that the mild mixing of operators lies strictly between strong mixing and weak mixing. 
By Proposition \ref{distinguish-examples} and Corollary \ref{hyper-not-dsj}, we immediately obtain that the mild mixing of invertible operators lies strictly between strong mixing and weak mixing. 

\begin{coro}\ 
\begin{enumerate}
    \item There exists an invertible $\mathcal{F}_{\mathrm{ld1}}$-transitive operator which is not mildly mixing. 
    \item There exists an invertible mildly mixing operator $T$ which is not strongly mixing. 
\end{enumerate}
\end{coro}

\begin{proof}
(1) The $\mathcal{F}_{\mathrm{ld1}}$-transitive operator $S\colon Y\to Y$ in Corollary \ref{hyper-not-dsj} is not mildly mixing since it is not weakly disjoint from some weakly mixing operator. 

(2) By Proposition \ref{distinguish-examples}(7), there exists an invertible $\Delta^*_{\mathrm{ip}}$-transitive operator $T\colon X\to X$ which is not $\mathcal{F}^*_{\mathrm{ip}}$-transitive. 
By Proposition \ref{mild-ip}, $T$ is mildly mixing. 
It is clear that $T$ is not strongly mixing. 
\end{proof}

We provide a further  example of a weakly mixing operator with a weakly mixing adjoint. 
By Remark \ref{TT^*}, this is also an example of a weakly mixing operator $T\colon X\to X$ such that neither $T$ nor $T^*$ is mildly mixing. 

\begin{coro}
There exists an invertible $\mathcal{F}_{\mathrm{ld1}}$-transitive operator $T$ on $\ell^2(\mathbb{Z})$ such that $T^*$ is weakly mixing. 
\end{coro}

\begin{proof}
    Since there exist $F_1\in\mathcal{F}_{\mathrm{ld1}}$ and $F_2\in\mathcal{F}_{\mathrm{t}}$ with $F_1\cap F_2=\emptyset$ (see the proof of Proposition \ref{distinguish-examples}(4)), the result follows from Proposition \ref{adjoint}. 
\end{proof}

\section{Direct sums of backward shifts}

In this section, we investigate the $\mathcal F$-transitivity of the direct sums of finitely many (bilateral and unilateral) backward shifts on Fr\'echet sequence spaces, 
and show that the $\mathcal{F}$-transitivity of such operators is equivalent to the $\mathcal{F}$-Mixing Criterion. 

Recall that in 1995 Salas gave in \cite{S95} a characterization of the weak disjointness of weighted backward shifts on $\ell^2(\mathbb{Z})$ and $\ell^2(\mathbb{N}_0)$. 
Indeed, the weak disjointness of weighted backward shifts on $c_0(\mathbb{Z})$ and $\ell^p(\mathbb{Z}), 1\leq p<\infty$ (or $c_0(\mathbb{N}_0)$ and $\ell^p(\mathbb{N}_0), 1\leq p<\infty$) is already characterized by Salas in \cite{S95}, 
although he only gave a proof for $\ell^2(\mathbb{Z})$ and $\ell^2(\mathbb{N}_0)$. 
We restate these results in the language of weak disjointness. 

\begin{thm}[{\cite[Theorem 2.5]{S95}}]\label{Salas-bilateral}
      Two weighted backward shifts $B_w, B_v$ on $\ell^2(\mathbb{Z})$ are weakly disjoint if and only if for any $\epsilon>0$ and $q\in\mathbb{N}_0$, there exists $n\in\mathbb{N}$ such that for any $|j|\leq q$, 
      \[
          \operatorname{min}\{w_{j+1}w_{j+2}\cdots w_{j+n}, v_{j+1}v_{j+2}\cdots v_{j+n}\}>\tfrac{1}{\epsilon},
      \]
      and
      \[
          \operatorname{max}\{w_{j}w_{j-1}\cdots w_{j-n+1}, v_{j}v_{j-1}\cdots v_{j-n+1}\}<\epsilon.
      \]
\end{thm}

\begin{thm}[{\cite[Theorem 2.8]{S95}}]\label{Salas}
    Two weighted backward shifts $B_w$ and $B_v$ on $\ell^2(\mathbb{N}_0)$ are weakly disjoint if and only if 
    \[
        \sup_{n\geq 1}\bigl\{\min\{\Pi_{s=1}^nw_s, \Pi_{s=1}^nv_s\}\bigr\}=\infty.
    \]
\end{thm}

In \cite{G00}, Grosse-Erdmann gave a characterization of weak disjointness of finitely many (unilateral and bilateral) backward shifts on F-sequence spaces with OP basis, by the limits of the basis, see \cite[Theorem 6, Theorem 7]{G00}. 

Let $X$ be a bilateral Fr\'echet sequence space in which $\{\dotsc, e_{-2},e_{-1}, e_0, e_1, e_2, \dotsc\}$ is a basis.
Given a bilateral weighted backward shift $B_w$ on $X$, it is not hard to see that $B_w$ is conjugate to a (non-weighted) backward shift $B$ on a Fr\'echet sequence space $X_w$ with the basis $\{\dotsc,w_{-1}w_0e_2,w_0e_{-1},e_0, \tfrac{e_1}{w_1}, \tfrac{e_2}{w_1w_2}, \dotsc\}$. 
This strategy also works for unilateral weighted backward shifts. 
Since the weak disjointness and the $\mathcal{F}$-transitivity of operators are preserved under conjugacy, in what follows we consider the direct sums of backward shifts, which can be transferred into the form of weighted backward shifts. 

\begin{thm}\label{F-dsj-r-bi}
    Let $X_1, X_2, \dotsc , X_k$ be bilateral Fr\'echet sequence spaces in which $(e_n)_{n\in\mathbb{Z}}$ is a basis and $B_1, B_2, \dotsc, B_k$ be the backward shifts on $X_1, X_2, \dotsc , X_k$ respectively. 
    Let $\mathcal{F}$ be a family. Then the following assertions are equivalent: 
     \begin{enumerate}
        \item $\bigoplus_{i=1}^kB_i$ is $\mathcal{F}$-transitive; 
        \item $\bigoplus_{i=1}^kB_i$ is $\mathcal{F}$-mixing; 
        \item $\bigoplus_{i=1}^kB_i$ satisfies the $\mathcal{F}$-Mixing Criterion; 
        \item for each $\epsilon>0$ and $N\in\mathbb{N}$, $\bigcap_{j=-N}^N\bigcap_{i=1}^k\bigl(C^+_{\epsilon, j}(X_i)\cap C^-_{\epsilon, j}(X_i)\bigr)\in\mathcal{F}$, where for $j\in\mathbb{Z}$ and $i\in\{1,2,\dotsc, k\}$, 
        \[
            C^+_{\epsilon, j}(X_i)=\bigl\{n\in\mathbb{N}_0: \|e_{j+n}\|_{X_i}<\epsilon\bigr\}, C^-_{\epsilon, j}(X_i)=\bigl\{n\in\mathbb{N}_0: \|e_{j-n}\|_{X_i}<\epsilon\bigr\}.
        \]
    \end{enumerate}
\end{thm}

\begin{proof}
The implications (3)$\Rightarrow$(2)$\Rightarrow$(1) are clear. 

(1)$\Rightarrow$(4). Fix $\epsilon>0$ and $N\in\mathbb{N}$. 
For each $i\in\{1,2,\dotsc,k\}$, since $(e_n)_{n\in\mathbb{Z}}$ is a basis of $X_i$, we have that $\{h_ne_n\colon n\in\mathbb{Z}\}$ is bounded for each $h=(h_n)_{n\in\mathbb{Z}}\in X_i$. 
By the Banach-Steinhaus Theorem (\cite[Theorem A.10]{G11}), for each $i\in\{1,2,\dotsc,k\}$, there exists $\delta_i>0$ such that for any $x^{(i)}\in X_i$ with $\|x^{(i)}\|_{X_i}<\delta_i$, we have that  
    \[
        \bigl\|x^{(i)}_ne_n\bigr\|_{X_i}<\frac{\epsilon}{2},\quad  n\in\mathbb{Z}, 
    \]
    and 
   \[
     \bigl|x^{(i)}_j\bigr|\leq\frac{1}{3},\quad |j|\leq N.
   \]
    Let $\delta=\operatorname{min}\{\delta_i\colon i=1,2,\dotsc k\}$ and $u=\sum_{i=-N}^Ne_i$. For each $i\in\{1,2,\dotsc,k\}$, denote $U_i=\{y\in X_i\colon \|y-u\|_{X_i}<\delta\}$. Since $\bigoplus_{i=1}^k B_i$ is hypercyclic, we have 
    \[
        Q:=\bigcap_{i=1}^k N_{B_i}(U_i,U_i)\cap [2N+1,+\infty)\neq\emptyset. 
    \]
    Fix $n\in Q$. There is a vector $(x^{(1)}, x^{(2)},\dotsc x^{(k)})\in U_1\oplus U_2\oplus\dotsc \oplus U_k$ such that 
    \[
       (B_1^nx^{(1)}, B_2^nx^{(2)},\dotsc B_k^nx^{(k)})\in U_1\oplus U_2\oplus\dotsc \oplus U_k. 
    \]
That is, we have $|j-n|>N, j+n>N$ for $|j|\leq N$ and for each $i\in\{1,2,\dotsc,k\}$, we have the following: 
    \begin{enumerate}[(i)]
        \item $|x^{(i)}_j-1|\leq\tfrac{1}{3}$ for $|j|\leq N$; 
        \item $|x^{(i)}_k|\leq\tfrac{1}{3}$ for $|k|>N$; 
        \item $\|(x^{(i)}_j-1)e_j\|_{X_i}<\tfrac{\epsilon}{2}$ for $|j|\leq N$; 
        \item $\|x^{(i)}_ke_k\|_{X_i}<\tfrac{\epsilon}{2}$ for $|k|>N$; 
        \item $|x^{(i)}_{n+j}-1|\leq\tfrac{1}{3}$ for $|j|\leq N$; 
        \item $|x^{(i)}_{n+k}|\leq\tfrac{1}{3}$ for $|k|>N$; 
        \item $\|(x^{(i)}_{n+j}-1)e_j\|_{X_i}<\tfrac{\epsilon}{2}$ for $|j|\leq N$; 
        \item $\|x^{(i)}_{n+k}e_k\|_{X_i}<\tfrac{\epsilon}{2}$ for $|k|>N$.
    \end{enumerate}
    
Fix $i\in\{1,2,\dotsc k\}$. By (viii), we have that $\|x^{(i)}_je_{j-n}\|_{X_i}<\frac{\epsilon}{2}$ for $|j|\leq N$. 
Then by (ii), we have that for $|j|\leq N$, 
\begin{align*}
    \|e_{j-n}\|_{X_i}&=\Bigl\|\bigl(\tfrac{1}{x^{(i)}_j}-1\bigr)x^{(i)}_je_{j-n}+ x^{(i)}_je_{j-n}\Bigr\|_{X_i}\\
    &\leq \bigr|\tfrac{1}{x^{(i)}_j}-1\bigr|\cdot  \bigl\|x^{(i)}_je_{j-n}\bigr\|_{X_i}+ \bigl\|x^{(i)}_je_{j-n}\bigr\|_{X_i}\\
    &<\epsilon.
\end{align*}
By (iv), we have that $\|x^{(i)}_{j+n}e_{j+n}\|_{X_i}<\frac{\epsilon}{2}$ for $|j|\leq N$. 
Then by (v), we have that for $|j|\leq N$, 
    \begin{align*}
        \|e_{j+n}\|_{X_i}&=\Bigl\|\bigl(\tfrac{1}{x^{(i)}_{j+n}}-1\bigr)x^{(i)}_{j+n}e_{j+n}+ x^{(i)}_{j+n}e_{j+n}\Bigr\|_{X_i} \\
        &\leq \Bigl\|\bigl(\tfrac{1}{x^{(i)}_{j+n}}-1\bigr)x^{(i)}_{j+n}e_{j+n}\Bigr\|_{X_i}+ \Bigl\|x^{(i)}_{j+n}e_{j+n}\Bigr\|_{X_i} \\
        &\leq \Bigl\|x^{(i)}_{j+n}e_{j+n}\Bigr\|_{X_i}+ \Bigl\|x^{(i)}_{j+n}e_{j+n}\Bigr\|_{X_i}\\
        &<\epsilon. 
    \end{align*}
Hence we have that
    \[
        Q \subset \bigcap_{j=-N}^N\bigcap_{i=1}^k\bigl(C^+_{\epsilon, j}(X_i)\cap C^-_{\epsilon, j}(X_i)\bigr).
    \]
Since $\bigoplus_{i=1}^kX_i$ has no isolated points, for each $i\in\{1,2,\dotsc,k\}$ there are non-empty open subsets $W_i$ and $V_i$ of $U_i$ such that
    \[
        \bigcap_{i=1}^k N_{B_i}(W_i, V_i)\subset  [2N+1, +\infty).
    \]
Then 
\[
    \bigcap_{i=1}^k N_{B_i}(W_i, V_i)\subset Q \subset \bigcap_{j=-N}^N\bigcap_{i=1}^k\bigl(C^+_{\epsilon, j}(X_i)\cap C^-_{\epsilon, j}(X_i)\bigr).
\]
We are done since the hypothesis gives that $\bigcap_{i=1}^k N_{B_i}(W_i, V_i)\in\mathcal{F}$. 

(4)$\Rightarrow$(3). Let $a=(a^{(1)}, a^{(2)}, \dotsc ,a^{(k)})$ and $b=(b^{(1)},b^{(2)}, \dotsc ,b^{(k)})$ be two $k$-tuples of finite sequences in $\bigoplus_{i=1}^k X_i$. 
For each $i\in \{1,2,\dotsc, k\}$, write 
\begin{align*}
    &a^{(i)}=(\dotsc, 0, 0, a^{(i)}_{-l}, \dotsc, a^{(i)}_0, a^{(i)}_1, \dotsc, a^{(i)}_l, 0, 0, \dotsc),\\
    &b^{(i)}=(\dotsc, 0, 0, b^{(i)}_{-l}, \dotsc, b^{(i)}_0, b^{(i)}_1, \dotsc, b^{(i)}_l, 0, 0, \dotsc).
\end{align*}
Let $F_1,F_2\dotsc, F_k$ be the forward shifts for finite sequences of $X_1,X_2\dotsc, X_k$ respectively. 
Denote $T=\bigoplus_{i=1}^k B_i$ and $S=\bigoplus_{i=1}^k F_i$.
Fix $\delta>0$. 
For each $i\in\{1,2,\dotsc,k\}$, there exists $\epsilon_i>0$ such that for any $z^{(i)}\in X_i$ with $\|z^{(i)}\|_{X_i}<\epsilon_i$, we have for $|j|\leq l$, 
\[
    \|a^{(i)}_jz^{(i)}\|_{X_i}<\frac{\delta}{2l+1}, \quad \|b^{(i)}_jz^{(i)}\|_{X_i}<\frac{\delta}{2l+1}. 
\]
Let $\epsilon=\operatorname{min}\{\epsilon_i\colon i=1,2,\dotsc,k\}$. 
Fix $n\in \bigcap_{j=-l}^l\bigcap_{i=1}^k\bigl(C^+_{\epsilon, j}(X_i)\cap C^-_{\epsilon, j}(X_i)\bigr)$ and $i\in\{1,2,\dotsc,k\}$. We have 

\[
    \bigl\|F_i^{n}a^{(i)}\bigr\|_{X_i}=\biggl\|\sum_{j=-l}^la^{(i)}_je_{j+n}\biggr\|_{X_i}\leq \sum_{j=-l}^l \bigl\|a^{(i)}_je_{j+n}\bigr\|_{X_i}<(2l+1) \frac{\delta}{2l+1}=\delta, 
\]
and 
\[
    \bigl\|B_i^{n}b^{(i)}\bigr\|_{X_i}=\biggl\|\sum_{j=-l}^lb^{(i)}_je_{j-n}\biggr\|_{X_i}\leq \sum_{j=-l}^l\bigl\|b^{(i)}_je_{j-n}\bigr\|_{X_i}<(2l+1) \frac{\delta}{2l+1}=\delta.
\]
Then we obtain: 
\[
    \mathcal{F}\text{-}\displaystyle\lim_{n\to\infty}(T^nb, S^na, T^nS^na)=(0,0,a).
\]
Hence $T$ satisfies the $\mathcal{F}$-Mixing Criterion. 
\end{proof}

We now turn to unilateral backward shifts. 
The corresponding result is stated below. 
Since its proof closely parallels that of Theorem \ref{F-dsj-r-bi} in the bilateral setting, we omit the details.

\begin{thm}\label{F-dsj-r-uni}
        Let $X_1, X_2, \dotsc, X_k$ be unilateral Fr\'echet sequence spaces in which $(e_n)_{n\in\mathbb{N}_0}$ is a basis and $B_1, B_2, \dotsc, B_k$ be the backward shifts on $X_1, X_2, \dotsc , X_k$ respectively. 
        Let $\mathcal{F}$ be a family. Then the following assertions are equivalent: 
     \begin{enumerate}
        \item $\bigoplus_{i=1}^kB_i$ is $\mathcal{F}$-transitive; 
        \item $\bigoplus_{i=1}^kB_i$ is $\mathcal{F}$-mixing; 
        \item $\bigoplus_{i=1}^kB_i$ satisfies the $\mathcal{F}$-Mixing Criterion; 
        \item for each $\epsilon>0$ and $N\in\mathbb{N}$, $\bigcap_{j=0}^N\bigcap_{i=1}^kC_{\epsilon, j}(X_i)\in\mathcal{F}$, where for $j\in\mathbb{N}_0$ and $i\in\{1,2,\dotsc,k\}$,
        \[
            C_{\epsilon, j}(X_i)=\bigl\{n\in\mathbb{N}_0: \|e_{n+j}\|_{X_i}<\epsilon\bigr\}.
        \]
    \end{enumerate}
\end{thm}

\noindent \textbf{Acknowledgments}: 
The authors would like to thank the referee for her or his valuable comments, which have helped to improve the presentation of this paper.
The authors were partially supported by the National Key R\&D Program of China (2024YFA1013601), NSF of China (12222110, 12171298) and a grant from the Guangdong Provincial Department of Education (2025KCXTD013).

\end{document}